\input amstex

\def\b1{\text{\bf 1}}

\def\BC{{\Bbb C}}
\def\BK{{\Bbb K}}

\def\BP{{\Bbb P}}
\def\cp{\BC\BP^{1}}
\def\cpn{\BC\BP^{N}}

\def\BZ{{\Bbb Z}}
\def\BV{{\Bbb V}}
\def\CA{{\Cal A}}

\def\CC{{\Cal C}}
\def\CD{{\Cal D}}

\def\CM{{\Cal M}}

\def\CO{{\Cal O}}
\def\CP{{\Cal P}}

\def\CT{{\Cal T}}
\def\CV{{\Cal V}}
\def\CW{{\Cal W}}

\def\dpar{\partial}
\def\End{\text{End}}

\def\fg{{\frak g}}

\def\Fields{\text{Fields}}

\def\hCT{\widehat{{\Cal T}}}
\def\hfg{\widehat{\frak g}}

\def\hsl{\widehat{\text{sl}}}
\def\htau{\widehat{\tau}}

\def\Id{\text{Id}}

\def\Lie{\text{Lie}}
\def\tlie{\widetilde{\Lie}}
\def\lie{{\Cal L}ie}
\def\Res{\text{Res}}
\def\Spec{\text{Spec}}

\def\ta{\tilde a}
\def\tb{\tilde b}

\def\tCT{\widetilde{\Cal{T}}}
\def\td{\tilde d}

\def\tf{\tilde f}

\def\tg{\tilde g}
\def\tI{\tilde I}

\def\tK{\tilde K}

\def\Vir{{\Cal V}ir}
\def\Mod{\text{Mod}}


\def\iso{\buildrel\sim\over\longrightarrow}

\def\lra{\longrightarrow}

\parskip=6pt

\documentstyle{amsppt}
\document
\NoBlackBoxes


\centerline{\bf Chiral de Rham complex. II}

\bigskip
\centerline{Fyodor Malikov, Vadim Schechtman}

\bigskip
\bigskip

\bigskip

{\it To Dmitry Borisovich Fuchs, on his sixtieth birthday}

\bigskip

\bigskip

\bigskip

\bigskip

This note is a sequel to [MSV]. It consists of three parts.
The first part is an expanded version of the last section of {\it op. cit.}
We give here certain construction of vertex algebras which
includes in particular the ones appearing in the above note.

In the second part we show how the cohomology ring $H^*(X)$
of a smooth complex variety $X$ could be restored from the
correlation functions of the vertex algebra $R\Gamma(X;\Omega^{ch}_X)$.

In the third part, we prove first a useful
general statement that the sheaf of loop algebras
over the tangent sheaf $\CT_X$ acts naturally on $\Omega^{ch}_X$
for every smooth $X$ (see \S 1).
The $\BZ$-graded vertex algebra $H^*(X;\Omega^{ch}_X)$
seems to be a quite interesting object (especially for compact $X$). 
In \S 2, we compute
$H^0(\cpn;\Omega^{ch}_{\cpn})$ as a module over $\hsl(N+1)$.

We are grateful to A. Beilinson and B. Feigin for useful
discussions.

\newpage

\bigskip
\centerline{\bf Part I: Chiral Weyl modules}
\bigskip

\bigskip

\bigskip
\centerline{\bf \S 1. Recollections on vertex algebras}
\bigskip

We will use the language of Kac's book [K] and of the original Borcherds'
paper [B]. All omitted proofs may be found in [K].

{\bf 1.1.} Let $V=V^{ev}\oplus V^{odd}$ be a super vector space. The
parity of an element $a\in V$ will be denoted by $\ta\in\BZ/2\BZ$.

We will denote by $V[[z,z^{-1}]]$
the space of all formal sums $f(z)=\sum_{i\in\BZ}\ a_iz^i;\ a_i\in V$, in
the even variable $z$. Similarly, we denote $V[[z,z^{-1},w,w^{-1}]]:=
V[[z,z^{-1}]][[w,w^{-1}]]$.

The subspace
of Laurent power series, with $a_i=0$ for $i<<0$ will be
denoted by $V((z))\subset V[[z,z^{-1}]]$ . We denote by
$\dpar_z:\ V[[z,z^{-1}]]\lra V[[z,z^{-1}]]$ the operator of differentiation
by $z$. We will also  write $f(z)'$ instead of $\dpar_zf(z)$.
If $A$ is a linear operator, $A^{(k)}$ will denote the operator $A^k/k!,\
k\in\BZ_{\geq 0}$.

{\bf 1.2.} A {\it vertex algebra} is a super vector space $V=V^{ev}\oplus
V^{odd}$ together with the following data.

(a) An element $|0\rangle\in V^{ev}$ called {\it vacuum vector}.

(b) An even linear map
$$
V\lra\End(V)[[z,z^{-1}]];\ a\mapsto a(z)=\sum_{n\in\BZ}\ a_{(n)}z^{-n-1}
\eqno{(1.2.1)}
$$
For each $a\in V$, the power series $a(z)$ must be a {\it field},
which means that for each $b\in V$, $a_{(n)}b=0$ for $n>>0$.
Here $a_{(n)}b$ denotes $a_{(n)}(b)\in V$. The coefficients $a_{(n)}$
are called {\it Fourier modes} of the field $a(z)$.

(c) An even endomorphism $\dpar:\ V\lra V$.

The following axioms must hold.

({\it Vacuum}) $|0\rangle(z)=\Id$ (the constant power series). $\dpar|0\rangle=
0$. For each $a\in V$, $a_{(n)}|0\rangle=0$ for $n\geq 0$, and
$a_{(-1)}|0\rangle=a$.

({\it Translation Invariance}) For each $a\in V$, $[\dpar,a(z)]=\dpar_za(z)$.

({\it Locality}) For all $a, b\in V$, there exists $N\in\BZ_{\geq 0}$
such that $(z-w)^N[a(z),b(w)]=0$.

Here
$$
[a(z),b(w)]:=\sum_{n,m}\ [a_{(n)},b_{(m)}]z^{-n-1}w^{-m-1}\in
\End(V)[[z,z^{-1},w,w^{-1}]]
$$

{\bf 1.3.} The following basic {\it Borcherds identity} follows from and together
with the Vacuum axiom is equivalent to the axioms of vertex algebra:
$$
\sum_{j=0}^\infty\ \binom{m}{j}\ (a_{(n+j)}b)_{(m+k-j)}=
$$
$$
=\sum_{j=0}^\infty\ (-1)^j\binom{n}{j}a_{(n+m-j)}b_{(k+j)}-
(-1)^{\ta\tb}\sum_{j=0}^\infty\ (-1)^{j+n}\binom{n}{j}b_{(n+k-j)}
a_{(m+j)}
\eqno{(1.3.1)}
$$
for each $a,b\in V,\ m,n,k\in\BZ$. Here by definition,
$$
\binom{n}{j}=\frac{n(n-1)\cdot\ldots\cdot(n-j+1)}{j!}
\eqno{(1.3.2)}
$$
for $n\in\BC, j\in\BZ_{\geq 0}$. Note that
$$
\binom{n}{j}=0\text{\ if\ }n\in\BZ, j>n\geq 0
\eqno{(1.3.3)}
$$
We set
$$
\binom{n}{j}=0\text{\ if\ }n,j\in\BZ; j<0
$$

Note the following useful particular case of (1.3.1) (corresponding to
$m=0$) which is in fact equivalent to (1.3.1),
$$
(a_{(n)}b)_{(k)}=\sum_{j=0}^\infty\ (-1)^j\binom{n}{j}a_{(n-j)}b_{(k+j)}-
(-1)^{\ta\tb}\sum_{j=0}^\infty\ (-1)^{j+n}\binom{n}{j}b_{(n+k-j)}a_{(j)}
\eqno{(1.3.4)}
$$

It is instructive to think of a vertex algebra as a super vector space
equipped with an infinite number of (nonassociative, noncommutative)  "multiplications"
$$
_{(n)}:\ a,b\mapsto a_{(n)}b\ \ (n\in\BZ)
\eqno{(1.3.5)}
$$
satisfying the quadratic relations (1.3.1).

Another important {\it commutativity formula} is
$$
a_{(n)}b=(-1)^{\ta\tb}\sum_{j=0}^\infty\ (-1)^{j+n+1}
\dpar^{(j)}(b_{(n+j)}a)
\eqno{(1.3.6)}
$$
for each $n\in\BZ$. Setting $n=0$ in (1.3.1), we get
$$
[a_{(m)},b_{(k)}]=\sum_{j=0}^\infty\ \binom{m}{j}(a_{(j)}b)_{(m+k-j)}
\eqno{(1.3.7)}
$$
for each $m,k\in\BZ$ ({\it Borcherds formula}). In particular,
$$
[a_{(0)},b_{(k)}]=(a_{(0)}b)_{(k)},
\eqno{(1.3.8)}
$$
in other words,
$$
a_{(0)}(b_{(k)}c)=(a_{(0)}b)_{(k)}c+(-1)^{\ta\tb}b_{(k)}(a_{(0)}c),
\eqno{(1.3.9)}
$$
i.e. the operators $a_{(0)}$ are derivations (of parity $\ta$) of all multiplications
$b_{(k)}c$.

It follows from the Vacuum and Translation invariance axioms that
$$
\dpar a=a_{(-2)}|0\rangle
\eqno{(1.3.10)}
$$
Substituting $m=0, n=-2, b=|0\rangle$ in (1.3.1) and taking into
account that
$$
\binom{-2}{j}=(-1)^j(j+1)
\eqno{(1.3.11)}
$$
we get
$$
(\dpar a)_{(k)}=-ka_{(k-1)},
\eqno{(1.3.12)}
$$
in other words,
$$
(\dpar a)(z)=\dpar_za(z)
\eqno{(1.3.13)}
$$
Iterating (1.3.9) we get
$$
a_{(-1-j)}=(\dpar^{(j)}a)_{(-1)}
\eqno{(1.3.14)}
$$
for all $j\in\BZ_{\geq 0}$. Applying (1.3.1) with $m=0, k=-2$ to the vacuum vector,
and taking into account (1.3.12), we deduce
$$
\dpar(a_{(n)}b)=(\dpar a)_{(n)}b+a_{(n)}(\dpar b),
\eqno{(1.3.15)}
$$
i.e. $\dpar$ is an even derivation of each multiplication $_{(n)}$.

{\bf 1.4. Language of fields.} Let $V$ be a vertex algebra. Let $\Fields(V)$
denote the subspace of $\End(V)[[z,z^{-1}]]$ consisting of all power
series $f(z)=\sum_{n\in\BZ}\ f_{(n)}z^{-n-1}$ such that for each $a\in V$
there exists $N\in\BZ$ such that $f_{(n)}(a)=0$ for each $n\geq N$.
This space inherits the obvious $\BZ/(2)$-grading from $\End(V)$.

We define a map $\pi:\ \Fields(V)\lra V$ by
$$
\pi(f(z))=f_{(-1)}(|0\rangle)
\eqno{(1.4.1)}
$$
By the vacuum axiom
$$
\pi(a(z))=a\ (a\in V)
\eqno{(1.4.2)}
$$
We set
$$
f(z)_-:=\sum_{n\geq 0}\ f_{(n)}z^{-n-1};\ f(z)_+=\sum_{n<0}\ f_{(n)}z^{-n-1}
\eqno{(1.4.3)}
$$
For two fields $f(z), g(z)$, we set
$$
:f(z)g(w):=f(z)_+g(w)+(-1)^{\tf\tg}g(w)f(z)_-\in\End(V)[[z,z^{-1},w,w^{-1}]]
\eqno{(1.4.4)}
$$
We can set $z=w$ in this expression and get the well-defined element
$:f(z)g(z):$ of $\Fields(V)$.
One defines the fields $f(z)_{(n)}g(z)\ (n\in\BZ)$ by the formulas
$$
f(z)_{(-1)}g(z)=:f(z)g(z):
\eqno{(1.4.5)}
$$
$$
f(z)_{(-1-j)}g(z)=\dpar_z^{(j)}(f(z)_{(-1)}g(z))\ \ (j\geq 0)
\eqno{(1.4.6)}
$$
(cf. (1.3.11)); and
$$
f(z)_{(j)}g(z)=\Res_w[f(w),g(z)](w-z)^j\ \ (j\geq 0)
\eqno{(1.4.7)}
$$
where $\Res_w$ denotes the coefficient at $w^{-1}$.

We have
$$
(a_{(n)}b)(z)=a(z)_{(n)}b(z)
\eqno{(1.4.8)}
$$
for each $a,b\in V,\ n\in\BZ$. Therefore,
$$
a_{(n)}b=\pi(a(z)_{(n)}b(z))
\eqno{(1.4.9)}
$$
The Borcherds formula (1.3.4) may be rewritten as
$$
a(z)b(w)=\sum_{j=0}^\infty\ \frac{(a_{(j)}b)(w)}{(z-w)^{j+1}}+
:a(z)b(w):
\eqno{(1.4.10)}
$$
This is understood as an identity in $V[[z,z^{-1},w,w^{-1}]]$; one
understands the fractions $(z-w)^{-j-1}$ as the elements of this space
using the binomial formula in the region $|z|>|w|$:
$$
\frac{1}{(z-w)^{j+1}}=\sum_{m=0}^\infty\ (-1)^m\binom{-j-1}{m}w^mz^{-j-1-m}=
$$
$$
=\sum_{n=0}^\infty\ \binom{n}{j}w^{n-j}z^{-n-1},
\eqno{(1.4.11)}
$$
cf. (1.3.3).
The identity (1.4.10) will be written as
$$
a(z)b(w)\sim\sum_{j=0}^\infty\ \frac{(a_{(j)}b)(w)}{(z-w)^{j+1}}
\eqno{(1.4.12)}
$$
(the {\it operator product expansion} of the fields $a(z), b(z)$).

{\bf 1.5.} A {\it morphism} of vertex algebras $f:\ V\lra V'$ is a
linear operator sending vacuum vector to the vacuum vector, and such that
$f(a_{(n)}b)=f(a)_{(n)}f(b)$ for all $a,b\in V,\ n\in\BZ$.

A linear operator $d:V\lra V$ of parity $\td$ is called a {\it derivation}
of the vertex algebra $V$ if
$$
d(a_{(n)}b)=d(a)_{(n)}b+(-1)^{\td\ta}a_{(n)}d(b)
\eqno{(1.5.1)}
$$
Thus, for each $a\in V$, $a_{(0)}$ is a derivation of $V$, and
$\dpar$ is an even derivation of $V$.

The {\it tensor product} $V\otimes W$ of two vertex algebras is their
tensor product as vector spaces with the vacuum vector
$|0\rangle_V\otimes |0\rangle_W$ and the following state-field correspondence
$$
(a\otimes b)(z)=a(z)\otimes b(z)
\eqno{(1.5.2)}
$$
i.e.
$$
(a\otimes b)_{(n)}=\sum_{k\in\BZ}\ a_{(k)}\otimes b_{(n-k-1)}
\eqno{(1.5.3)}
$$
It follows that
$$
\dpar_{V\otimes W}=\dpar_V\otimes\Id_W+\Id_V\otimes\dpar_W
\eqno{(1.5.4)}
$$

{\bf 1.6.} A {\it graded} vertex algebra is a vertex algebra $V$, together with
an even diagonalizable linear operator $H:\ V\lra V$ ({\it Hamiltonian}) such that
$$
[H,a(z)]=z\dpar_za(z)+(Ha)(z)
\eqno{(1.6.1)}
$$
for each $a\in V$. We will denote by $V_{\Delta}$ the eigenspace
of $H$ corresponding to the eigenvalue $\Delta$. The eigenvalues of $H$ are called
{\it conformal weights}. For $a\in V_\Delta$, we will write the field
$a(z)$ in the form
$$
a(z)=\sum_{i\in -\Delta+\BZ}\ a_iz^{-i-\Delta}
\eqno{(1.6.2)}
$$
Thus,
$$
a_i=a_{(i+\Delta-1)};\ a_{(n)}=a_{n-\Delta+1}
\eqno{(1.6.3)}
$$
The identity (1.6.1) means that $a_n$ has conformal weight $-n$, i.e.
$$
a_n(V_\Delta)\subset V_{\Delta-n}
\eqno{(1.6.4)}
$$
In other words,
$$
_{(n)}:\ V_{\Delta}\otimes V_{\Delta'}\lra V_{\Delta+\Delta'-n-1}
\eqno{(1.6.5)}
$$
In particular,
$$
_{(-1)}:\ V_{\Delta}\otimes V_{\Delta'}\lra V_{\Delta+\Delta'}
\eqno{(1.6.6)}
$$
It follows that
$$
|0\rangle\in V_0
\eqno{(1.6.7)}
$$
since $a_{(-1)}|0\rangle=a$, and
$$
\dpar(V_\Delta)\subset V_{\Delta+1}
\eqno{(1.6.8)}
$$
since $\dpar a=a_{(-2)}|0\rangle$.

The OPE formula (1.4.12) will be rewritten as
$$
a(z)b(w)\sim\sum_{j=0}^\infty\ \frac{(a_{j-\Delta+1}b)(w)}{(z-w)^{j+1}}\ \
(a\in V_\Delta)
\eqno{(1.6.9)}
$$
or
$$
[a_m,b(z)]=\sum_{j=0}^\infty\ \binom{m+\Delta-1}{j}(a_{j-\Delta+1}b)(z)\cdot
z^{m-j+\Delta-1}\ (a\in V_\Delta)
\eqno{(1.6.10)}
$$

{\bf 1.7.} A {\it conformal vertex algebra of central charge $c\in \BC$} is a vertex algebra $V$
equipped
with an even vector $L\in V$ satisfying the conditions (a) --- (c) below.

(a) The field $L(z)$ satisfies the OPE
$$
L(z)L(w)\sim\frac{c}{2(z-w)^4}+\frac{2L(w)}{(z-w)^2}+\frac{L(w)'}{z-w}
\eqno{(1.7.1)}
$$
Equivalently, if we write $L(z)=\sum\ L_nz^{-n-2}$ (so that $L_{(n)}=L_{n-1}$),
then the components $L_n$ satisfy the Virasoro commuation relations
$$
[L_m,L_n]=(m-n)L_{m+n}+\frac{m^3-m}{12}\cdot c\cdot \delta_{m,-n}
\eqno{(1.7.2)}
$$

(b) $L_{-1}=\dpar$

(c) The operator $L_0$ is diagonalizable.

It follows that $V$ is a graded vertex algebra with Hamiltonian $L_0$.

We have
$$
L(z)b(w)\sim\sum_{n=-1}^\infty\ \frac{(L_nb)(w)}{(z-w)^{n+2}}=
$$
$$
=\frac{b(w)'}{z-w}+\frac{\Delta b(w)}{(z-w)^2}+\ldots\ \ (b\in V_{\Delta})
\eqno{(1.7.3)}
$$
or
$$
[L_m,b(z)]=\sum_{j=0}^\infty\ \binom{m+1}{j}(L_{j-1}b)(z)\cdot z^{m-j+1}=
$$
$$
=z^{m+1}b(z)'+(m+1)\Delta z^mb(z)+\ldots\ \ (b\in V_\Delta)
\eqno{(1.7.4)}
$$

{\bf 1.8.} A {\it conformal algebra} is a superspace $R=R^{ev}\oplus R^{odd}$,
together with an even linear operator $\dpar:\ R\lra R$ and a collection
of even operations $_{(n)}:\ R\otimes R\lra R,\ n\in\BZ_{\geq 0}$
satisfying the axioms (C0) --- (C3) below.

(C0) For each $a,b\in R$, $a_{(n)}b=0$ for $n>>0$.

(C1) $(\dpar a)_{(n)}b=-na_{(n-1)}b$

(C2) $a_{(n)}b=(-1)^{\ta\tb}\sum_{j=0}^\infty\ \dpar^{(j)}(b_{(n+j)}a)$

(C3) $a_{(m)}(a_{(n)}c)=\sum_{j=0}^\infty\ \binom{m}{j}
(a_{(j)}b)_{(m+n-j)}c+(-1)^{\ta\tb}b_{(n)}(a_{(m)}c)$

One has an obvious forgetful functor
$$
(Vertex\ algebras)\lra (Conformal\ algebras)
\eqno{(1.8.1)}
$$
This functor admits a left adjoint
$$
U:\ (Conformal\ algebras)\lra (Vertex\ algebras)
\eqno{(1.8.2)}
$$
called the {\it vertex envelope}, cf. [K], 4.7.

In the language of [BD], vertex algebras correspond to {\it (unital) chiral algebras}
over the disk (we have the equivalence of categories).
Conformal algebras correspond to {\it $Lie^*$-algebras}. The analogue
of the functor $U$ is called the {\it chiral envelope}.

{\bf 1.9.} A vertex algebra $V$ is called {\it holomorphic} if
$a_{(n)}=0$ for each $a\in V,\ n\geq 0$.

If $V$ is a holomorphic vertex algebra, then the operation $ab:=a_{(-1)}b$
is supercommutative and associative, and the vacuum vector is a unity. The operator
$\dpar$ is an even derivation, $\dpar(ab)=(\dpar a)b+a\dpar b$. The remaining
operations $_{(-n-1)},\ n\geq 0,$ can be recovered by the formula
(1.3.14).

This
gives an equivalence of the categories of holomorphic vertex algebras and
the category of commutative associative unital superalgebras with an
even derivation, cf. [K], 1.4.

In the language of [BD], holomorphic vertex algebras correspond to
{\it commutative} chiral algebras over the disk.
The above mentioned equivalence is
translated into the equivalence of the category of commutative chiral
algebras and that of $\CD$-algebras (commutative
algebras in the category of $\CD$-modules over the disk).

{\bf 1.10.} Let $V$ be a vertex algebra. A {\it module} over
$V$ is a super vector space $M$, together with an even linear map
$$
V\lra\ \End(M)[[z,z^{-1}]],\ \ a\mapsto a(z)=\sum_{n\in\BZ}\
a_{(n)}z^{-n-1}
\eqno{(1.10.1)}
$$
such that for all $a\in V,\ m\in M$, $a_{(n)}m=0$ for $n>>0$,
$|0\rangle(z)=\Id_M$, and the Borcherds identity (1.3.1) holds true.

If $V$ is graded, then $M$ is called graded if it is equipped with a
direct sum decomposition $M=\oplus\ M_{\Delta}$ such that

for all $a\in V_{\Delta},\ m\in M_{\Delta'},\ n\in\BZ$, we have
$a_{(n)}m\in M_{\Delta+\Delta'-n-1}$,

cf (1.6.5).

\bigskip
\centerline{\bf \S 2. Restricted vertex algebras}
\bigskip

{\bf 2.1.} Let $(V,H)$ be a graded vertex algebra. Consider the condition

(P) $V$ has no negative integer conformal weights.

Vertex algebras satisfying this condition will be called {\it restricted}.

{\bf 2.2.} Let us fix a restricted vertex algebra $V$. To simplify
the formulas, we will assume that $V$ is purely even. All the considerations
below have the obvious "super" ($\BZ/(2)$-graded) version.

The operation $a_{(-1)}b$ will be denoted simply by $ab$ and referred to as a multiplication by
$a$.

We are going to investigate, what kind of a structure
on the subspace $V_{\leq 1}:=V_0\oplus V_1$
is induced by our vertex algebra.

(a) The space $V_0$ is a commutative assotiative unital $\BC$-algebra
with respect to the operation $ab$. The unity equals $|0\rangle$.

Note that
$$
a_{(n)}b=0\text{\ for\ } a,b\in V_0,\ n\geq 0
\eqno{(2.2.1)}
$$
since the operation $_{(n)}$ has conformal weight $-n-1$, cf. (1.6.10).
Now, the commutativity of $V_0$ follows from the commutativity formula
(1.3.6).

Applying (1.3.1) with $m=0,\ n=k=-1$, to an element $c$, and taking into
account that
$$
\binom{-1}{j}=(-1)^j
\eqno{(2.2.2)}
$$
we get
$$
(ab)c=a(bc)+\sum_{j=0}^\infty\ a_{(-2-j)}b_{(j)}c+
\sum_{j=0}^\infty\ b_{(-2-j)}a_{(j)}c
\eqno{(2.2.3)}
$$
for all $a, b, c$. If $a, b, c\in V_0$ then the sums disappear,
and we get the associativity.
The Vacuum axiom implies that $|0\rangle$ is a unity.

This algebra will be denoted $A$. The vacuum will be denoted $\b1$.

There is a map $A\otimes V\lra V,\ a\otimes b\mapsto ab$.
However,
it does not make $V$ into an $A$-module: the multiplication by $A$ is
not associative in general.

We have the map $\dpar:\ A\lra V_1$. Let $\Omega\subset V_1$
denote the subspace spanned by the elements $a\dpar b,\ a,b\in A$.
Thus, $\dpar$ induces the map
$$
d:\ A\lra\Omega
\eqno{(2.2.4)}
$$

(b) The left multiplication by $A$ makes $\Omega$ a left $A$-module.
We have
$$
a db=(db)a\ \ (a,b\in A)
\eqno{(2.2.5)}
$$

Let us write down a particular case of (2.2.3):
$$
(ab)c-a(bc)=a_{(-2)}b_{(0)}c+b_{(-2)}a_{(0)}c
$$
$$
\text{\ for\ }a,b\in V_0,\ c\in V_1
\eqno{(2.2.6)}
$$
For $b,c\in V_0$, we have $0=\dpar(b_{(0)}c)=(\dpar b)_{(0)}c+b_{(0)}\dpar c$,
but $(\dpar b)_{(0)}=0$ by (1.3.12). It follows that
$$
b_{(0)}\dpar c=0\text{\ for\ }b,c\in V_0
\eqno{(2.2.7)}
$$
Therefore,
$$
(ab)\dpar c=a(b\dpar c)\text{\ for\ }a,b,c\in V_0
\eqno{(2.2.8)}
$$
It follows that $\Omega$ is a left $A$-module.

(c) The map $d$ is a derivation, i.e.
$$
d(ab)=adb+bda
\eqno{(2.2.9)}
$$

This follows from (1.3.15).

Let us denote by $\CT$ the quotient space $V_1/\Omega$.

(d) The left multiplication by $A$ makes $\CT$ into a left $A$-module.

Note that by (1.3.12),
$$
a_{(-2)}b=\dpar a\cdot b
\eqno{(2.2.10)}
$$
Thus, we can rewrite (2.2.3) as
$$
(ab)c-a(bc)=\dpar a\cdot (b_{(0)}c)+\dpar b\cdot (a_{(0)}c)
\eqno{(2.2.11)}
$$

Consider the operation
$$
_{(0)}:\ V_1\otimes V_1\lra V_1
\eqno{(2.2.12)}
$$

(e) The operation (2.2.12) induces a Lie bracket on $\CT$, to be denoted
$[\ ,\ ]$.

By
(1.3.6), we have
$$
a_{(0)}b=-b_{(0)}a+\dpar(b_{(1)}a)
\eqno{(2.2.13)}
$$
Therefore, composition
$$
V_1\otimes V_1\buildrel{(0)}\over\lra V_1\lra \CT
\eqno{(2.2.14)}
$$
is skew-symmetric. On the other hand, by (1.3.1) with $m=0, n=-1, k=0$,
we have
$$
(a\dpar b)_{(0)}c=-\dpar a\cdot (b_{(0)}c)+\dpar b\cdot (a_{(0)}c)
\text{\ for\ }a, b\in V_0,\ c\in V_1
\eqno{(2.2.15)}
$$
It follows that the composition (2.2.14) is zero on the subspace
$\Omega\otimes V_1+V_1\otimes \Omega$. The Jacobi identity
holds since the operators $a_{(0)}$ are derivations of all
operations $_{(k)}$.

We will omit the proof of claims (f) --- (n) below; they are proved in a similar
manner by application of Borcherds identity and commutativity
formula. The reader may want to perform these (easy) calculations on its own.

Consider the operation
$$
_{(0)}:\ V_1\otimes A\lra A
\eqno{(2.2.16)}
$$

(f) The operation (2.2.16) vanishes on the subspace $\Omega\otimes A$
and induces on $A$ a structure of a module over the Lie algebra $\CT$.

This action will be denoted by $\tau(a)\ (a\in A,\ \tau\in \CT)$.

(g) The Lie algebra $\CT$ acts on $A$ by derivations,
$$
\tau(ab)=\tau(a)b+a\tau(b)\
\eqno{(2.2.17)}
$$

(h) We have
$$
(a\tau)(b)=a\tau(b)
\eqno{(2.2.18)}
$$

The properties (d) --- (h) mean that $\CT$ is a {\it Lie algebroid} over $A$,
cf. [BFM], 3.2.1.

(i) The operation (2.2.12) induces on the space $\Omega$ a structure of a module over the
Lie algebra $\CT$.

This action will be denoted $\tau(\omega)$ or
$\tau\omega\ (\tau\in \CT,\ \omega\in\Omega)$.

(j) We have
$$
\tau(a\omega)=a\tau(\omega)+\tau(a)\omega\ (a\in A,\ \tau\in\CT,\
\omega\in \Omega)
\eqno{(2.2.19)}
$$

(k) The differential $d:\ A\lra\Omega$ is compatible with the
$\CT$-module structure.

It follows from (j) and (k) that

(l) we have
$$
\tau(a db)=\tau(a)db + a d(\tau(b))\ \ (\tau\in\CT,\ a,b\ \in A)
\eqno{(2.2.20)}
$$

Consider the operation
$$
_{(1)}:\ V_1\otimes V_1\lra A
\eqno{(2.2.21)}
$$

(m) The map (2.2.21) vanishes on the subspace
$\Omega\otimes\Omega$. Therefore, it induces the pairing
$$
\langle\ ,\ \rangle:\ \Omega\otimes\CT\oplus\CT\otimes\Omega\lra A
\eqno{(2.2.22)}
$$
This pairing is $A$-bilinear and symmetric.
We have
$$
\langle \tau, a db\rangle=a\tau(b)\ \ (\tau\in\CT,\ a,b\in A)
\eqno{(2.2.23)}
$$

(n) We have
$$
(a\tau)(\omega)=a\tau(\omega)+\langle\tau,\omega\rangle da\ \ (a\in A,\ \tau\in\CT,\
\omega\in\Omega)
\eqno{(2.2.24)}
$$

{\bf 2.3.} Let us denote by $\hCT$ the space $V_1/dA$.
The operation (2.2.4) induces a Lie bracket on the space $\hCT$.
The subspace $\Omega/dA\subset\hCT$ is an abelian Lie ideal.
The adjoint action of $\CT=\hCT/(\Omega/dA)$ coincides with
the action defined by (i) and (l).

Thus, we have an extension of Lie algebras
$$
0\lra \Omega/dA\lra \hCT\lra \CT\lra 0
\eqno{(2.3.1)}
$$
Note that this extension is not central in general. 

Let us denote the space $V_1$ by $\tCT$. We can form a dg Lie algebra 
$$
A\lra \tCT\lra \hCT
\eqno{(2.3.2)}
$$
living in degrees $-2, -1, 0$. The first arrow is $d$, the second 
one is the projection. The $(-1,-1)$-component of the bracket is given 
by the operation (2.2.21), the $(0,-1)$-component is induced 
by the operation $_{(0)}$. 

{\bf 2.3.1. Remark.} In fact, the dg Lie algebra (2.3.2) is a part 
of a bigger dg Lie algebra which one can associate with an 
arbitrary vertex algebra $V$. It is defined as 
$$
V\lra V\lra V/\dpar V
\eqno{(2.3.1.1)}
$$ 
and lives in degrees $-2, -1, 0$. The first arrow is $\dpar$. 
The operation $_{(0)}$ induces the $V/\dpar V$-module structure on $V$; 
this will be the $(0,-1)$-component of the bracket. The $(-1,-1)$-component 
is given by the symmetric operation  
$$
[a,b]^{-1,-1}=\sum_{j\geq 0}\ (-1)^j\frac{\dpar^j}{(j+1)!}a_{(j+1)}b
\eqno{(2.3.1.2)}
$$
cf. [B], Section 9. The $(0,-2)$-component is trivial.  
This definition was inspired by [BFM], 6.4.    

{\bf 2.4.}
We have an exact sequence of vector spaces
$$
0\lra\Omega\lra\tCT\lra\CT\lra 0
\eqno{(2.4.1)}
$$
Both arrows are compatible with the left multiplication by $A$. Let $\pi$
denote the projection $\pi:\ \tCT\lra\CT$.

Let us define the "bracket" $[\ ,\ ]:\ \Lambda^2\tCT\lra\tCT$ by
$$
[x,y]=\frac{1}{2}(x_{(0)}y-y_{(0)}x)\ \ (x,y\in\tCT)
\eqno{(2.4.2)}
$$
This bracket does not make $\tCT$ into a Lie algebra: the Jacobi identity
is in general violated. Set
$$
J(x,y,z)=[[x,y],z]+[[y,z],x]+[[z,x],y]\ \ (x,y,z\in \tCT)
\eqno{(2.4.3)}
$$
Consider the operation (2.2.21).

(a) The operation (2.2.21) is symmetric. It will be denoted
by $\langle x,y\rangle$.

Let us define the map $I:\ \Lambda^3\tCT\lra A$ by
$$
I(x,y,z)=\langle x,[y,z]\rangle+\langle y,[z,x]\rangle+
\langle z,[x,y]\rangle
\eqno{(2.4.4)}
$$

(b) We have
$$
J(x,y,z)=\frac{1}{6}dI(x,y,z)
\eqno{(2.4.5)}
$$

In fact, by (1.3.6),
$$
[x,y]=x_{(0)}y-\frac{1}{2}d\langle x,y\rangle
\eqno{(2.4.6)}
$$
It follows that
$$
[[x,y],z]=[x,y]_{(0)}z-\frac{1}{2}d\langle [x,y],z\rangle=
$$
$$
=(x_{(0)}y)_{(0)}z-\frac{1}{2}d\langle [x,y],z\rangle=
x_{(0)}y_{(0)}z-y_{(0)}x_{(0)}z-\frac{1}{2}d\langle [x,y],z\rangle
$$
Similarly,
$$
[[y,z],x]=(y_{(0)}z)_{(0)}x-\frac{1}{2}d\langle [y,z],x\rangle,
$$
and a similar identity for $[[z,x],y]$. Now we use the identities
$$
(y_{(0)}z)_{(0)}x=-x_{(0)}y_{(0)}z+d\langle x,y_{(0)}z\rangle,
$$
and the similar one for $(z_{(0)}x)_{(0)}y$, to get
$$
J(x,y,z)=d\bigl(\langle x,y_{(0)}z\rangle -\langle y,x_{(0)}z\rangle\bigr)-
\frac{1}{2}d\bigl(\langle x,[y,z]\rangle + (cycle)\bigr)
$$
Note that the left hand side of this equality, as well as the second
summand in the right hand side (equal to $-1/2\cdot dI(x,y,z)$), are
completely skew symmetric
with respect to all permutations of the letters $x,y,z$. Therefore,
we can skew symmetrize the first summand as well,
which gives $2/3\cdot dI(x,y,z)$. The identity (2.4.5) follows.

(c) We have
$$
\langle ax,y\rangle=a\langle x,y\rangle-\pi(x)\pi(y)(a)\ \ (a\in A,\
x,y\in\tCT)
\eqno{(2.4.7)}
$$

This follows from (1.3.1)  with $m=0, n=-1, k=1$.

(d) We have
$$
\langle [x,y],z\rangle+\langle y,[x,z]\rangle =\pi(x)(\langle y,z\rangle)-
\frac{1}{2}\pi(y)(\langle x,z\rangle) -
\frac{1}{2}\pi(z)(\langle x,y\rangle)
\eqno{(2.4.8)}
$$

To check this, we use the identity
$$
x_{(0)}\langle y,z\rangle=\langle x_{(0)}y,z\rangle +
\langle y,x_{(0)}z\rangle
\eqno{(2.4.9)}
$$
and take into account that
$$
\langle da,z\rangle=\pi(z)(a)\ \ (a\in A,\ z\in \tCT)
\eqno{(2.4.10)}
$$

{\bf 2.5.} Let us choose a splitting
$$
s:\ \CT\lra\tCT
\eqno{(2.5.1)}
$$
of the projection $\pi$. Let us define the map
$$
\langle\ ,\ \rangle=\langle\ ,\ \rangle_s:\ S^2\CT\lra A
\eqno{(2.5.2)}
$$
by
$$
\langle\tau_1,\tau_2\rangle=\langle s(\tau_1),s(\tau_2)\rangle
\eqno{(2.5.3)}
$$
(we use the lower index $_s$ in the notation if we want to stress
the dependence on the splitting $s$).
Let us define the map
$$
c^2=c^2_s:\ \Lambda^2\CT\lra \Omega
\eqno{(2.5.4)}
$$
by
$$
c^2(\tau_1,\tau_2)=s([\tau_1,\tau_2])-[s(\tau_1),s(\tau_2)]
\eqno{(2.5.5)}
$$
Let us define the map $K:\ \Lambda^3\CT\lra A$ by
$$
K(\tau_1,\tau_2,\tau_3)=\langle s(\tau_1),s([\tau_2,\tau_3])\rangle+
\langle s(\tau_2),s([\tau_3,\tau_1])\rangle+
\langle s(\tau_3),s([\tau_1,\tau_2])\rangle\
\eqno{(2.5.6)}
$$
Let us define the map
$$
c^3=c^3_s:\ \Lambda^3\CT\lra A
$$
by
$$
c^3(\tau_1,\tau_2,\tau_3)=-\frac{1}{2}K(\tau_1,\tau_2,\tau_3)
+\frac{1}{3}I(s(\tau_1),s(\tau_2),s(\tau_3))
\eqno{(2.5.7)}
$$
cf. (2.4.4). Let us regard $c^2$ (resp. $c^3$) as a Lie algebra cochains
living in $C^2(\CT;\Omega)$ (resp., in $C^3(\CT;A)$).
Recall that the Lie differential $d_{Lie}$ is defined by
$$
d_{Lie}(c^i)(\tau_1,\ldots,\tau_{i+1})=
\sum_{p<q}\ (-1)^{p+q-1}c^i([\tau_p,\tau_q],\tau_1,\ldots,\htau_p,\ldots,
\htau_q,\ldots,\tau_{i+1})+
$$
$$
+\sum_p\ (-1)^p\tau_pc^i(\tau_1,\ldots,\htau_p,\ldots,\tau_{i+1})
\eqno{(2.5.8)}
$$

(a) We have
$$
d_{Lie}(c^2)=dc^3
\eqno{(2.5.9)}
$$

In fact, one sees from the definition (2.5.5) that
$$
d_{Lie}(c^2)(\tau_1,\tau_2,\tau_3)=-[s([\tau_1,\tau_2]),s(\tau_3)]+(cycle)
-s(\tau_1)_{(0)}s([\tau_2,\tau_3])+(cycle)+
$$
$$
+s(\tau_1)_{(0)}[s(\tau_2),s(\tau_3)]+(cycle) -
J(s(\tau_1),s(\tau_2),s(\tau_3))
$$
Now, (2.5.9) follows from (2.4.6) and (2.4.5).

(b) We have
$$
d_{Lie}(c^3)=0
\eqno{(2.5.10)}
$$

To check this, note that $c^3$ is a sum of two summands, as in (2.5.7).
Also, let us split the map $d_{Lie}$ into $d'+d''$ where
$$
d'(\gamma)(\tau_1,\ldots)=\sum_{p<q}\ (-1)^{p+q-1}
\gamma([\tau_p,\tau_q],\tau_1,\ldots,\htau_p,\ldots,\htau_q,
\ldots)
$$
and
$$
d''(\gamma)(\tau_1,\ldots)=\sum_p\ (-1)^p\tau_p \gamma(\tau_1,\ldots,
\htau_p,\ldots)
$$
Correspondingly, $d_{Lie}(c^3)$ splits into four summands.

The summand $-\frac{1}{2}d'K$ contains terms of the type
$$
\langle s([\tau_1,\tau_2]),s[\tau_3,\tau_4])\rangle
\eqno{(2.5.11)}
$$
(six terms) and the terms of the type
$$
\langle s(\tau_1),s([\tau_2,[\tau_3,\tau_4]])\rangle,
$$
these terms cancel out, due to the Jacobi identity in $\CT$.

The summand $-\frac{1}{2}d''K$ has the terms of the type
$$
-\tau_1(\langle s(\tau_2),s([\tau_3,\tau_4])\rangle+
\tau_2(\langle s(\tau_1),s([\tau_3,\tau_4])
\eqno{(2.5.12)}
$$
(six groups). The following identity is an easy consequence of (2.4.8):
$$
\frac{3}{2}\bigl\{ \pi(x)(\langle y,z\rangle)-\pi(y)(\langle x,z\rangle)\bigr\}=
2\langle [x,y],z\rangle + \langle y,[x,z]\rangle -
\langle x,[y,z]\rangle
\eqno{(2.5.13)}
$$
Therefore, (2.5.12) is equal to
$$
\frac{2}{3}\bigl\{-2\langle [s(\tau_1),s(\tau_2)],s([\tau_3,\tau_4])\rangle-
\langle s(\tau_2),[s(\tau_1),s([\tau_3,\tau_4])]\rangle+
$$
$$
+\langle s(\tau_1),[s(\tau_2),s([\tau_3,\tau_4])]\rangle\bigr\}
\eqno{(2.5.14)}
$$

The summand $\frac{1}{3}d'I$ contains six groups of the type
$$
\langle s([\tau_3,\tau_4]),[s(\tau_1),s(\tau_2)]\rangle+
\langle s(\tau_1),[s(\tau_2),s([\tau_3,\tau_4])]\rangle+
$$
$$
+\langle s(\tau_2),[s([\tau_3,\tau_4]),s(\tau_1)]\rangle
\eqno{(2.5.15)}
$$
The second the the third terms in (2.5.15) cancel out with the similar terms
in (2.5.14).

Using (2.5.13) and (2.4.5), it is easy to deduce the identity
$$
\tau_1(\langle s(\tau_2),[s(\tau_3),s(\tau_4)]\rangle)+\ldots=
3\bigl\{
\langle [s(\tau_1),s(\tau_2)],[s(\tau_3),s(\tau_4)]\rangle-
$$
$$
-\langle [s(\tau_1),s(\tau_3)],[s(\tau_2),s(\tau_4)]\rangle+
\langle [s(\tau_1),s(\tau_4)],[s(\tau_2),s(\tau_3)]\rangle\bigr\}
\eqno{(2.5.16)}
$$
where $\ldots$ in the left hand side mean the complete skew-symmetrization
(twelve summands altogether).
It follows that
$$
\frac{1}{3}d''I=
-\langle [s(\tau_1),s(\tau_2)],[s(\tau_3),s(\tau_4)]\rangle+
\langle [s(\tau_1),s(\tau_3)],[s(\tau_2),s(\tau_4)]\rangle-
$$
$$
-\langle [s(\tau_1),s(\tau_4)],[s(\tau_2),s(\tau_3)]\rangle
\eqno{(2.5.17)}
$$
Finally, since
$$
\langle [s(\tau_1),s(\tau_2)],c(\tau_3,\tau_4)\rangle=
\langle s([\tau_1,\tau_2]),c(\tau_3,\tau_4)\rangle,
$$
we have
$$
\langle s([\tau_1,\tau_2]),s([\tau_3,\tau_4])\rangle=
\langle [s(\tau_1),s(\tau_2)],s([\tau_3,\tau_4])\rangle+
\langle s([\tau_1,\tau_2]),[s(\tau_3),s(\tau_4)]\rangle-
$$
$$
-\langle [s(\tau_1),s(\tau_2)],[s(\tau_3),s(\tau_4)]\rangle
\eqno{(2.5.18)}
$$
Using this relation, we see that the sum of the remaining summands equals
zero. This proves (2.5.10).

The identities (a) and (b) mean that the pair $c=(c^2,c^3)$ form a
$2$-cocycle of the Lie algebra $\CT$ with coefficients in the
complex $A\lra\Omega$.

(c) We have
$$
\langle[\tau_1,\tau_2],\tau_3\rangle+\langle\tau_2,[\tau_1,\tau_3]\rangle=
\tau_1(\langle\tau_2,\tau_3\rangle)-
\frac{1}{2}\tau_2(\langle\tau_1,\tau_3\rangle)-
\frac{1}{2}\tau_3(\langle\tau_1,\tau_2\rangle)-
$$
$$
-\langle \tau_2,c^2(\tau_1,\tau_3)\rangle -
\langle\tau_3,c^2(\tau_1,\tau_2)\rangle
\eqno{(2.5.19)}
$$
This follows from 2.4 (d).

Let us investigate the effect of the change of a splitting.
Let $s':\ \CT\lra \tCT$ be another splitting of $\pi$. The difference
$s'-s$ lands in $\Omega$; let us denote it
$$
\omega=\omega_{s,s'}:\ \CT\lra\Omega
\eqno{(2.5.20)}
$$
We regard $\omega$ as a $1$-cochain of $\CT$ with coefficients in $\Omega$.
Let us define a $2$-cochain $\alpha=\alpha_{s,s'}\in C^2(\CT;A)$ by
$$
\alpha(\tau_1,\tau_2)=\frac{1}{2}\bigl(\langle \omega(\tau_1),\tau_2\rangle-
\langle \tau_1,\omega(\tau_2)\rangle\bigr)
\eqno{(2.5.21)}
$$

(d) We have
$$
c^2_s-c^2_{s'}=d_{Lie}(\omega)-d\alpha
\eqno{(2.5.22)}
$$

Indeed, the left hand side of (2.5.22) is equal to
$$
\omega([\tau_1,\tau_2])-[s(\tau_1),\omega(\tau_2)]-[\omega(\tau_1),
s'(\tau_2)]
$$
By (2.4.6), we have
$$
[s(\tau_1),\omega(\tau_2)]=\tau_1\omega(\tau_2)-\frac{1}{2}d(\langle \tau_1,
\omega(\tau_2)\rangle),
$$
and the similar expression for $[\omega(\tau_1),s'(\tau_2)]$.
The identity (2.5.22) follows.

(e) We have
$$
-d_{Lie}(\alpha)=c^3_s-c^3_{s'}
\eqno{(2.5.23)}
$$

Indeed we have
$$
d_{Lie}\alpha(\tau_1,\tau_2,\tau_3)=
\frac{1}{2}\bigl\{\langle\omega([\tau_1,\tau_2]),\tau_3\rangle -
\ldots -\langle [\tau_2,\tau_3],\omega(\tau_1)\rangle\bigr\}+
$$
$$
+\frac{1}{2}\bigl\{-\tau_1(\langle\omega(\tau_2),\tau_3\rangle)+\ldots
+\tau_3(\langle\tau_1,\omega(\tau_2)\rangle)\bigr\}
$$
and
$$
c^3_s(\tau_1,\tau_2,\tau_3)-c^3_{s'}(\tau_1,\tau_2,\tau_3)=
-\frac{1}{2}\bigl\{\langle\omega(\tau_1),[\tau_2,\tau_3]\rangle+
\ldots + \langle\tau_3,\omega([\tau_1,\tau_2])\rangle\bigr\}+
$$
$$
+\frac{1}{3}\bigl\{\langle\omega(\tau_1),[\tau_2,\tau_3]\rangle+
\ldots +\langle\tau_3,[\tau_1,\omega(\tau_2)]\rangle\bigr\}
\eqno{(2.5.24)}
$$
The summands of the type $\langle\tau_1,\omega([\tau_2,\tau_3])\rangle$
in (2.5.24) are equal to the corresponding summands in the expression
for $-d_{Lie}\alpha$. It remains to consider the summands containing
$\omega(\tau_i)$. Consider for example the summands in (2.5.24)
containing $\omega(\tau_1)$. There is one,
$$
-\frac{1}{2}\langle\omega(\tau_1),[\tau_2,\tau_3]\rangle
$$
in the first bracket, and three,
$$
\frac{1}{3}\bigl\{\langle\omega(\tau_1),[\tau_2,\tau_3]\rangle+
\langle\tau_2,[\omega(\tau_3),\tau_1]\rangle+
\langle\tau_3,[\tau_1,\omega(\tau_2)]\rangle\bigr\}
$$
in the second one. In the above sum, replace the second term using
(2.4.8),
$$
\langle\tau_3,[\omega(\tau_1),\tau_2]\rangle=
\langle [\tau_2,\tau_3],\omega(\tau_1)]\rangle -
\tau_2(\langle\tau_3,\omega(\tau_1)\rangle)+
\frac{1}{2}\tau_3(\langle\tau_2,\omega(\tau_1)\rangle),
$$
and similarly replace the third term $\langle\tau_3,[\tau_1,\omega(\tau_2)]
\rangle$. Summing up, we get the same as in $-d_{Lie}\alpha$.

The properties (d) and (e) mean that
$$
c_s-c_{s'}=d_{Lie}(\beta)
\eqno{(2.5.25)}
$$
where $\beta=\beta_{s,s'}:=(\omega,\alpha)\in C^1(\CT;A\lra\Omega)$.

Therefore, we have assigned to our vertex algebra a canonically defined
"characteristic class"
$$
c(V)=c(V_{\leq 1})\in H^2(\CT;A\lra\Omega)
\eqno{(2.5.26)}
$$
as the cohomology class of the cocycle $c_s$.

{\bf 2.6.} Note that we have
$$
(ab)x-a(bx)=-\pi(x)(a)db - \pi(x)(b)da\ \ (a,b\in A,\ x\in\tCT)
\eqno{(2.6.1)}
$$
This follows from (1.3.1) with $m=0, n=-1, k=-1$.

Let us introduce the mapping
$$
\gamma=\gamma_s:\ A\otimes\CT\lra\Omega
\eqno{(2.6.2)}
$$
so that
$$
\gamma(a,\tau)=s(a\tau)-a s(\tau)
\eqno{(2.6.3)}
$$

(a) We have
$$
\gamma(ab,\tau)=\gamma(a,b\tau)+a\gamma(b,\tau)+
\tau(a)db+\tau(b)da
\eqno{(2.6.4)}
$$

This follows from (2.6.1).

We have
$$
(ax)_{(0)}y=a(x_{(0)}y)-\pi(y)(a)x+
\langle x,y\rangle da-d(\pi(x)\pi(y)(a))
\eqno{(2.6.5)}
$$
$(a\in A,\ x,y\in\tCT)$, which follows from (1.3.1) with $m=0, n=-1, k=0$.
This implies (using (2.4.6) and (2.4.7))
$$
[ax,y]=a[x,y]-\pi(y)(a)x+\frac{1}{2}\langle x,y\rangle da -
\frac{1}{2}d(\pi(x)\pi(y)(a))
\eqno{(2.6.6)}
$$
$(a\in A,\ x,y\in\tCT)$.

(b) We have
$$
\langle a\tau_1,\tau_2\rangle=a\langle\tau_1,\tau_2\rangle +
\langle \gamma(a,\tau_1),\tau_2\rangle -\tau_1\tau_2(a)
\eqno{(2.6.7)}
$$

This follows from (2.4.7).

(c) We have
$$
c^2(a\tau_1,\tau_2)=ac^2(\tau_1,\tau_2)+\gamma(a,[\tau_1,\tau_2])-
\gamma(\tau_2(a),\tau_1)+\tau_2\gamma(a,\tau_1)-
$$
$$
-\frac{1}{2}\langle\tau_1,\tau_2\rangle da+
\frac{1}{2}d(\tau_1\tau_2(a))-
\frac{1}{2}d(\langle\tau_2,\gamma(a,\tau_1)\rangle)
\eqno{(2.6.8)}
$$
$(a\in A,\ \tau_i\in\CT)$.

This follows from (2.6.6) and (2.6.7).

(d) We have
$$
c^3(a\tau_1,\tau_2,\tau_3)=ac^3(\tau_1,\tau_2,\tau_3)+
\frac{1}{2}\tau_1[\tau_2,\tau_3](a)-
$$
$$
-\frac{1}{2}\bigl\{\langle\tau_2,\gamma(a,[\tau_3,\tau_1])\rangle-
\langle\tau_3,\gamma(a,[\tau_2,\tau_1])\rangle+
\langle\tau_2,\gamma(\tau_3(a),\tau_1)\rangle-
\langle\tau_3,\gamma(\tau_2(a),\tau_1)\rangle\bigr\}+
$$
$$
+\frac{1}{2}\bigl\{\langle\tau_2,\tau_3\gamma(a,\tau_1)\rangle-
\langle\tau_3,\tau_2\gamma(a,\tau_1)\rangle\bigr\}
-\frac{1}{2}\langle [\tau_2,\tau_3],\gamma(a,\tau_1)\rangle
\eqno{(2.6.9)}
$$

(e) {\bf Exercise.} Check that the formulas (c) and (d) are compatible
with the identity $d_{Lie}(c^2)=dc^3$.

Hint:  use 2.2 (n).

{\bf 2.7. Filtered algebras.} In fact, we do not 
really need the gradation on our vertex algebras. Let us call 
a vertex algebra $V$ {\it filtered} if it is eqipped with an 
increasing exhaustive filtration 
$\ldots\subset V_{\leq i}\subset V_{\leq i+1}\subset\ldots\subset V$ such that 
$|0\rangle\in V_{\leq 0}$ and  
$$
V_{\leq i}\ _{(n)}V_{\leq j}\subset V_{\leq i+j-n-1}
\eqno{(2.7.1)}
$$
cf. (1.6.5). As a consequence, $\dpar V_{\leq i}\subset V_{\leq i+1}$. 
(The corresponding graded space $gr V$ is a graded vertex algebra.)   

Let us call $V$ {\it restricted} if $V_{\leq -1}=0$.  
For such $V$, again $A:=V_{\leq 0}$ is a commutative algebra. 
Define $\Omega\subset V_{\leq 1}$ as the $\BC$-subspace spanned 
by $a\dpar b,\ a,b\in A$. Set 
$$
\CA:=V_{\leq 1}/\Omega;\ \CT:=V_{\leq 1}/(A+\Omega)
\eqno{(2.7.2)}
$$
Then $\CA$ is an $A$-Lie algebroid ({\it Atiyah algebra} of $V$) which 
is an $A$-extension 
of the Lie algebroid $\CT$, cf. [BFM].  
The discussion 2.2 --- 2.6 carries over to the filtered case, 
with $\CT$ replaced by $\CA$.     

\bigskip

\newpage

\bigskip
\centerline{\bf \S 3. Prevertex algebras}
\bigskip

{\bf 3.1.} Let us call the data (a) --- (f) below a {\bf prevertex algebra}.

(a) A commutative algebra $A$.

(b) An $A$-module $\Omega$, together with an $A$-derivation
$d:\ A\lra \Omega$. We assume that $\Omega$ is generated as a vector
space by the elements $a db\ (a,b\in A)$, i.e. the canonical map
$\Omega^1(A):=\Omega^1_{\BC}(A)\lra\Omega$ is surjective.

(c) An $A$-Lie algebroid $\CT$. Define the action of $\CT$ on $\Omega$ by
$$
\tau(a db)=\tau(a) db+ a d(\tau(b)),
\eqno{(3.1.1)}
$$
cf. (2.2.9). We assume that this action is well defined.
It follows that $d$ is compatible with the action of $\CT$.

We assume that the formula
$$
\langle\tau,a db \rangle=a\tau(b)
\eqno{(3.1.2)}
$$
gives a well defined $A$-bilinear pairing $\CT\times\Omega\lra A$.

(d) A $\BC$-bilinear mapping $\gamma:\ A\times\CT\lra\Omega$ satisfying
2.6 (a).

(e) A $\BC$-bilinear symmetric mapping $\langle\ ,\ \rangle:\
\CT\times\CT\lra A$ satisfying 2.6 (b).

(f) A skew symmetric map $c^2:\ \Lambda^2\CT\lra \Omega$. This map should
satisfy 2.5 (c), 2.6 (c), and the property (3.1.7) below. Let us define the map
$$
[\ ,\ ]':= [\ ,\ ]-c^2:\ \Lambda^2\CT\lra \CT\oplus\Omega
\eqno{(3.1.3)}
$$
Let us define the map
$$
c^3:=-\frac{1}{2}\tK+\frac{1}{3}\tI:\ \Lambda^3\CT\lra A
\eqno{(3.1.4)}
$$
where
$$
\tK(\tau_1,\tau_2,\tau_3)=\langle\tau_1,[\tau_2,\tau_3]\rangle+
\langle\tau_2,[\tau_3,\tau_1]\rangle+
\langle\tau_3,[\tau_1,\tau_2]\rangle
\eqno{(3.1.5)}
$$
and
$$
\tI(\tau_1,\tau_2,\tau_3)=\langle\tau_1,[\tau_2,\tau_3]'\rangle+
\langle\tau_2,[\tau_3,\tau_1]'\rangle+
\langle\tau_3,[\tau_1,\tau_2]'\rangle
\eqno{(3.1.6)}
$$
In the last formula, we extend the symmetric pairing $\langle\ ,\ \rangle$
to the whole space $\CT\oplus\Omega$, by using (3.1.2), (e), and
setting it equal to zero on $\Omega\times\Omega$.

Now, with $c^3$ defined above, we require that
$$
d_{Lie}(c^2)=dc^3;\ d_{Lie}(c^3)=0
\eqno{(3.1.7)}
$$

{\bf 3.2.}  Given a prevertex algebra $P=(A,\Omega,\CT,\ldots)$, set $V_0=A;\
V_1=\CT\oplus\Omega$. Define the operation
$_{(-2)}:\ V_0\times V_0\lra V_1$ by
$$
a_{(-2)}b=b da
\eqno{(3.2.1)}
$$

Define the operations
$_{(-1)}:\ V_0\times V_1\lra V_1$ and $_{(-1)}:\ V_1\times V_0\lra V_1$ by
$$
a_{(-1)}\omega=a\omega;\
a_{(-1)}\tau=a\tau+\gamma(a,\tau)
\eqno{(3.2.2)}
$$
and
$$
\omega_{(-1)}a=a\omega;\ \tau_{(-1)}a=a\tau+\gamma(a,\tau)+
d(\tau(a))
\eqno{(3.2.3)}
$$
$(a\in A, \omega\in\Omega, \tau\in\CT)$.

Define the operations $_{(0)}:\ V_0\times V_1\lra V_0$ and
$_{(0)}:\ V_1\times V_0\lra V_0$ by
$$
\tau_{(0)}a=\tau(a),\
a_{(0)}\tau=-\tau(a),\ \omega_{(0)}a=a_{(0)}\omega=0
\eqno{(3.2.4)}
$$

Define the operation
$_{(0)}:\ V_1\times V_1\lra V_1$ by
$$
\omega_1\ _{(0)}\omega_2=0;\
\tau_{(0)} a db=\tau(a) db + a d(\tau(b));\
(a db)_{(0)}\tau=-\tau(a) db + \tau(b) da
\eqno{(3.2.5)}
$$
and
$$
\tau_1\
_{(0)}\tau_2=[\tau_1,\tau_2]-c^2(\tau_1,\tau_2)+\frac{1}{2}d\langle\tau_1,\tau_2\rangle
\eqno{(3.2.6)}
$$

Define the symmetric operation $_{(1)}:\ V_1\times V_1\lra V_0$ by
$$
\omega_1\ _{(1)}\omega_2=0;\
\tau_{(1)}(a db)=a\tau(b);\
\tau_1\ _{(1)}\tau_2=\langle\tau_1,\tau_2\rangle
\eqno{(3.2.7)}
$$

This defines a structure of a {\it one-truncated vertex algebra}
on the space $V_{\leq 1}:=V_0\oplus V_1$. This means that
we have the vacuum vector in $V_0$,
the map $\dpar:=d:\ V_0\lra V_1$, together with the operations
$$
_{(n)}:\ V_i\times V_j\lra V_{i+j-n-1}
$$
for all $i,j,n$ such that $i,j,i+j-n-1\in\{0,1\}$. These operations
satisfy the Borcherds identities whenever the indices belong to
the indicated range, the impossible operations being set to zero.

Let us give examples of prevertex algebras.

{\bf 3.3. Example.} Let $\CT$ be a Lie algebra over $\BC$ equipped
with a symmetric invariant bilinear form $\langle\ ,\ \rangle$.
Set $A=\BC$, $\Omega=0,\ c^2=0, \gamma=0$. This defines a prevertex algebra.

Note that the component $c^3$ defined by the rule (f) is equal to
$$
c^3(\tau_1,\tau_2,\tau_3)=-\frac{1}{2}\langle\tau_1,[\tau_2,\tau_3]\rangle
\eqno{(3.3.1)}
$$

{\bf 3.4. Example.} Let $A$ be a $\BC$-algebra; set
$\Omega=\Omega^1_{\BC}(A)$. Let $\CT_0$ be an abelian Lie algebra over
$\BC$ acting by derivations on $A$.

(For example, let $A$ be smooth, and assume that there exists a basis
$\{\tau_i\}$
of the left $A$-module $\CT:=Der_{\BC}(A)$ consisting of commuting
vector fields. Let $\CT_0$ be the vector space spanned by this basis.)

Set $\CT=A\otimes_{\BC}\CT_0$. There is a unique Lie bracket on $\CT$
making it a Lie algebroid over $A$,
$$
[a\tau_1,b\tau_2]=a\tau_1(b)\tau_2-b\tau_2(a)\tau_1\ \
(\tau_i\in\CT_0, a, b\in A)
\eqno{(3.4.1)}
$$

We set $\langle\tau_1,\tau_2\rangle=0;\ \gamma(a,\tau)=0;
c^2(\tau_1,\tau_2)=0; c^3(\tau_1,\tau_2,\tau_3)=0$ for
$a\in A,\ \tau_i\in\CT_0$. Then the formulas 2.6 (a) --- (d) define a
unique extension of the operations $\gamma, \langle\ ,\ \rangle,
c^2$ and $c^3$ to the whole space $\CT$.

Namely,
$$
\gamma(a,b\tau)=-\tau(a)db-\tau(b)da
\eqno{(3.4.2)}
$$
$$
\langle a\tau_1,b\tau_2\rangle=-a\tau_2\tau_1(b)-b\tau_1\tau_2(a)-
\tau_1(b)\tau_2(a)
\eqno{(3.4.3)}
$$
It is convenient to write down $c=(c^2,c^3)$ as a sum of a simpler
cocycle and a coboundary,
$$
c^2(a\tau_1,b\tau_2)=\ 'c^2(a\tau_1,b\tau_2)+d\beta(a\tau_1,b\tau_2)
\eqno{(3.4.4)}
$$
where
$$
'c^2(a\tau_1,b\tau_2)=\frac{1}{2}\bigl\{\tau_1(b)d\tau_2(a)-
\tau_2(a)d\tau_1(b)\bigr\}
\eqno{(3.4.5)}
$$
$$
\beta(a\tau_1,b\tau_2)=
\frac{1}{2}\bigl\{b\tau_1\tau_2(a)-a\tau_2\tau_1(b)\bigr\}
\eqno{(3.4.6)}
$$
and
$$
c^3(a\tau_1,b\tau_2,c\tau_3)=\ 'c^3(a\tau_1,b\tau_2,c\tau_3)+
d_{Lie}\beta(a\tau_1,b\tau_2,c\tau_3)
\eqno{(3.4.7)}
$$
where
$$
'c^3(a\tau_1,b\tau_2,c\tau_3)=\frac{1}{2}\bigl\{
\tau_1(b)\tau_2(c)\tau_3(a)-\tau_1(c)\tau_2(a)\tau_3(b)\}
\eqno{(3.4.8)}
$$
A straightforward check shows that this gives a prevertex algebra.

Note that the cocycle $('c^2,'c^3)$ coincides with (minus one half of)
the one from [MSV], (5.17-18) if $A$ is the polynomial ring.

{\bf 3.5. Example.} Here is a common generalization of the
above two examples. Namely, in the set up of 3.4, assume that $\CT_0$ is an
arbitrary Lie algebra over $\BC$ acting by derivations on $A$, and
equipped with a symmetric invariant bilinear form $\langle\ ,\ \rangle$. There is a unique
extension of the Lie structure on $\CT_0$ to a Lie algebroid structure on $\CT:=
A\otimes_{\BC}\CT_0$, given by
$$
[a\tau_1,b\tau_2]=ab[\tau_1,\tau_2]+a\tau_1(b)\tau_2-b\tau_2(a)\tau_1
\eqno{(3.5.1)}
$$
Set $\gamma(a,\tau)=0$ for $\tau\in\CT_0, a\in A$. There exists
a unique extension of the function $\gamma$ to the whole $\CT$, given by
$$
\gamma(a,b\tau)=-\tau(a)db-\tau(b)da
\eqno{(3.5.2)}
$$
There exists a unique extension of the pairing $\langle\ ,\ \rangle$ to the whole
$\CT$ satisfying (2.6.7), namely
$$
\langle a\tau_1,b\tau_2\rangle=ab\langle\tau_1,\tau_2\rangle-a\tau_2\tau_1(b)-
b\tau_1\tau_2(a)-\tau_1(b)\tau_2(a)
\eqno{(3.5.3)}
$$
The cochain $c^2$ is defined as the unique extension of the zero cocycle on $\CT_0$
satisfying (2.6.8). It is given by
$$
c^2=\ 'c^2+\ ''c^2+d\beta
\eqno{(3.5.4)}
$$
where $\beta$ is as in (3.3.6),
$$
'c^2(a\tau_1,b\tau_2)=\frac{1}{2}\bigl\{\tau_1(b)d\tau_2(a)-\tau_2(a)d\tau_1(b)\bigr\}-
$$
$$
-[\tau_1,\tau_2](a)db-[\tau_1,\tau_2](b)da
\eqno{(3.5.5)}
$$
and
$$
''c^2(a\tau_1,b\tau_2)=\frac{1}{2}\langle\tau_1,\tau_2\rangle(a db-b da)
\eqno{(3.5.6)}
$$
(this is a "Kac-Moody" type summand). The component $c^3$ is given by
$$
c^3=\ 'c^3+\ ''c^3
\eqno{(3.5.7)}
$$
where
$$
'c^3(a\tau_1,b\tau_2,c\tau_3)=\frac{1}{2}\bigl\{\tau_1(b)\tau_2(c)\tau_3(a)-
\tau_1(c)\tau_2(a)\tau_3(b)\bigr\}
\eqno{(3.5.8)}
$$
and
$$
''c^3(a\tau_1,b\tau_2,c\tau_3)=-\frac{1}{2}abc\langle\tau_1,[\tau_2,\tau_3]\rangle
\eqno{(3.5.9)}
$$
cf. (3.3.1). We have
$$
d_{Lie}(\ 'c^2)=d(\ 'c^3);\ d_{Lie}(\ ''c^2)=d(\ ''c^3)
\eqno{(3.5.10)}
$$
and both $'c^3,\ ''c^3$ are Lie $3$-cocycles.

{\bf 3.6. Base change.} The previous example is a particular case of the
following construction. Let $P=(A,\Omega,\CT,\ldots)$ be a prevertex algebra,
and let $A'$ be a commutative $A$-algebra.

Assume that we are given the data (a), (b) below.

(a) An action of the Lie algebra $\CT$ on the algebra $A'$ by
derivations, extending its action on $A$ and such that
$$
(a\tau)(b)=a\tau(b)\ (a\in A,\ b\in A')
\eqno{(3.6.1)}
$$

(b) A mapping $\gamma:\ A'\times\CT\lra\Omega_{A'}$
extending the map 3.1 (d) such that 2.6 (a) is satisfied for all
$a\in A', b\in A$. Here $\Omega_{A'}:=A'\otimes_A\Omega$.

We call the data (a), (b) {\it the base change data}.

Set $\CT_{A'}:=A'\otimes_A\CT$. Due to (a), there is a unique
$A'$-Lie algebroid structure on $\CT_{A'}$, extending the given
$A$-Lie algebroid structure on $\CT$. Formulas (2.6.4), (2.6.7)
and (2.6.8) define the extention of the operations
$\gamma, \langle\ ,\ \rangle$ and $c^2$ to the space
$\CT_{A'}$.

This defines a new prevertex algebra
$P_{A'}=(A',\Omega_{A'},\CT_{A'},\ldots)$.

{\bf 3.7.} The additive group of the space
$Hom_{\BC}(\CT,\Omega)=C^1_{Lie}(\CT;\Omega)$ is acting on the
set of data 2.7 (a) --- (f) with fixed $A, \Omega, \CT, d$. Namely,
for $w:\ \CT\lra\Omega$, let us define the new
$$
\gamma_w(a,\tau)=\gamma(a,\tau)+w(a\tau)-aw(\tau)
\eqno{(3.7.1)}
$$
$$
\langle\tau_1,\tau_2\rangle_w=\langle\tau_1,\tau_2\rangle+
\langle w(\tau_1),\tau_2\rangle+\langle\tau_1,w(\tau_2)\rangle
\eqno{(3.7.2)}
$$
$$
c^2_w(\tau_1,\tau_2)=c^2(\tau_1,\tau_2)+d_{Lie}w(\tau_1,\tau_2)+
d\alpha(\tau_1,\tau_2)
\eqno{(3.7.3)}
$$
$$
c^3_w(\tau_1,\tau_2,\tau_3)=d_{Lie}\alpha(\tau_1,\tau_2,\tau_3)
\eqno{(3.7.4)}
$$
where
$$
\alpha(\tau_1,\tau_2)=\frac{1}{2}\bigl\{\langle\tau_1,w(\tau_2)\rangle-
\langle\tau_2,w(\tau_1)\rangle\bigr\}
\eqno{(3.7.5)}
$$
The resulting truncated vertex algebras are isomorphic;
if we denote the corresponding weight one spaces by
$V_1,\ V_{w1}$ where $V_{w1}=V_1$ as a vector space,
then the isomorphism $\phi_w:\ V_1\iso V_{w1}$ is defined by
$$
\phi_w(\tau,\omega)=(\tau,\omega+w(\tau))\ \ (\tau\in\CT,\ \omega\in\Omega)
\eqno{(3.7.6)}
$$
$\phi_w$ being the identity on $V_0=V_{0w}$. We have obviously
$$
\phi_{w+w'}=\phi_w\circ\phi_{w'}
\eqno{(3.7.7)}
$$
This means that the truncated vertex algebra is defined uniquely, up
to a unique isomorphism, by an element of the quotient of the set
of the above data by the action of $Hom_{\BC}(\CT,\Omega)$.

{\bf 3.8.} Let $V$ be a restricted vertex algebra, and let $\CM$ be
a graded $V$-module. We say that $\CM$ is restricted if it does not
contain the negative integer conformal weights.

Suppose that $\CM$ is restricted. Again we assume to simplify the notations
that $V$ and $M$ are even.
We will adopt the notations of 2.2.
We will denote the operation $a_{(-1)}m$ simply
by $am$. Consider the weight zero component $M=\CM_0$. The operation $am$
makes $M$ a left $A$-module.

The operation $_{(-1)}:\ V_1\otimes M\lra M$ vanishes on the subspace
$\Omega\otimes M$, and thus induces the operation
$$
\CT\otimes M\lra M,\ \ \tau\otimes m\mapsto \tau m
\eqno{(3.8.1)}
$$
making $M$ a module over the Lie algebra $\CT$.
This operation satisfies
$$
a(\tau m)=(a\tau) m
\eqno{(3.8.2)}
$$
and
$$
\tau(am)=a(\tau m)+\tau(a) m
\eqno{(3.8.3)}
$$
$(a\in A,\ \tau\in\CT,\ m\in M)$.

For example, if $\CT$ coincides with the algebra of vector fields,
this means that $M$ is a $\CD$-module over $A$.

\bigskip
\centerline{\bf \S 4. Vertex envelope}
\bigskip

{\bf 4.1.} We have a functor
$$
(Restricted\ vertex\ algebras)\lra (Commutative\ algebras)
\eqno{(4.1.1)}
$$
sending a vertex algebra $V$ to its conformal weight zero
component. To shorten the notation, by a "commutative algebra"
we mean a "commutative associative unital algebra".

We claim that this functor admits a left adjoint.

To construct it, we note that the forgetful functor from the
category of the commutative algebras with a
derivation to that of commutative algebras
admits a left adjoint. It sends an algebra $R$ to the algebra
$R_{\dpar}$ which is the quotient of the free commuative algebra over
$R$ generated by symbols $\dpar^{(n)}r,\ r\in R,\ n\in\BZ_{\geq 0}$,
over the relations
$$
\dpar^{(0)}r=r;\ \dpar^{(n)}(r_1+r_2)=\dpar^{(n)}r_1+\dpar^{(n)}r_2;
$$
$$
\dpar^{(n)}(r_1r_2)=\sum_{p+q=n}\dpar^{(p)}r_1\cdot\dpar^{(q)}r_2
\eqno{(4.1.2)}
$$
The derivation acts as $\dpar(\dpar^{(n)}r)=(n+1)\dpar^{(n+1)}r$.

Now, the left adjoint to (4.1.1) assigns to a commutative algebra $R$
the holomorphic vertex algebra corresponding to $R_\dpar$, cf. 1.8.

{\bf 4.2.} We have a functor
$$
(Restricted\ vertex\ algebras)\lra (One-truncated\ vertex\ algebras)
\eqno{(4.2.1)}
$$
sending a vertex algebra $V$ to $V_{\leq 1}=V_0\oplus V_1$.
This functor admits a left adjoint, $U_1$. Namely, 
given a one-truncated algebra $W$, we first form 
the non-commutative nonassociative algebra $W'$ spanned by 
$\dpar^{(n)}W,\ n\geq 0$ subject to relations (4.1.2). We 
introduce the operations $_{(n)}$ on $W'$ using (1.3.12) and 
(1.3.4). By definition, $U_1W$ is the quotient of $W'$ 
by the Borcherds relations (1.3.1). We have $(U_1W)_{\leq 1}=W$.
   
The size of $U_1W$ 
is that of the symmetric algebra on the vector space 
$$
\oplus_{n\geq 0}\dpar^nA\oplus\ \oplus_{n\geq 0}\dpar^n\CT
$$ 
(where $(A,\CT,\ldots)$ is the prevertex 
algebra corresponding to $W$). $U_1W$ contains 
as a subalgebra the holomorphic vertex 
algebra corresponding to $A$.  

{\bf Remark.} Essentially the same construction works in the filtered case, 2.7. 
Vertex algebras of this type appear as chiral algebras 
of {\it twisted} differential operators, cf. 4.4 below and 
[MSV], 5.15.  
  
Let us call a restricted vertex algebra $V$ {\it split} if it
is given together with a splitting (2.5.1). We have defined in
2.2 --- 2.6 a functor
$$
\CP:\ (Split\ restricted\ vertex\ algebras)\lra (Prevertex\ algebras)
\eqno{(4.2.2)}
$$
This functor admits a left adjoint, to be denoted by $\CV$,
and called {\it vertex envelope}. Namely, given 
a prevertex algebra, we first take the corresponding 
one-truncated vertex algebra, 3.2, and then apply $U_1$.

{\bf 4.3. Example.} If $P$ is as in 3.3 then $\CV(P)$ is
the vacuum (level $1$) representation of the affine Lie
algebra corresponding to $(\CT,\langle\ ,\ \rangle)$,
cf. [K], 4.7.

{\bf 4.4. Example.} If $P$ is as in 3.4,
we get the algebras studied in [MSV], cf. {\it op. cit.}, 6.9.

{\bf 4.5. Chiral Weyl modules.} Let us return again to the even situation.
Let $V$ be a restricted vertex algebra,
let $\CP(V)=(A,\Omega,\CT,\ldots)$ be the corresponding
prevertex algebra. Assume that the Lie
algebra $\CT$ coincides with $Der(A)$.

In 3.8 we have defined a functor
$$
(Restricted\ V-modules)\lra (D_A-modules)
\eqno{(4.5.1)}
$$
where $D_A$ is the algebra of differential operators on $A$.  
This functor admits a left adjoint
$$
\CW:\ (D_A-modules) \lra (Restricted\ V-modules)
\eqno{(4.5.2)}
$$
Its construction is similar to 4.2, and even simpler. 
We leave it to the reader. 

For a $D_A$-module $M$, the $V$-module $\CW(M)$ is called
the {\it chiral Weyl module} corresponding to $M$.

\newpage

\bigskip
 \centerline{\bf Part II: Vertex algebras and coinvariants }
\bigskip

\bigskip

In \S 1 we focus on two important features of an arbitrary conformal
vertex algebra $V$. Firstly, with any such an algebra and any
smooth curve we associate a sheaf of Lie algebras, see 1.1-1.3.
We use this sheaf to define a space of (co)invariants (or ``conformal
blocks'') in the situation when there are several $V$-modules attached
to several points on the curve; these spaces arrange in a vector bundle
with flat connection, see 1.6. In 1.8-9 we explain that
 one  can represent a horizontal section of the bundle associated with
several copies of $V$ attached to points on $\cp$ as a matrix element
of a product of fields.

Secondly, we show in 1.4
 that the Lie algebra of Fourier components of fields
associated with any conformal vertex algebra affords an antiinvolution.
This gives a duality functor on the category of $V$-modules, see 1.7.

In \S\S 2,3 we explain what all of this means in the case of the sheaf
$\Omega^{ch}_{X}$ introduced in [MSV]. In particular, we reproduce
the well-known physics calculation that represents the structure coefficients
of the cohomology ring of $X$ as matrix elements of products of fields,
see 3.2.

\bigskip

\centerline{\bf   \S 1. The spaces of (co)invariants}

\bigskip

{\bf 1.1.} We place ourselves in the situation of I.1.7, that is, assume that
we are given a vertex algebra $V$ along with a Virasoro field
$L(z)\in\Fields(V)$ satisfying (I.1.7.1-2). Consider the space of
Fourier modes of fields $\Lie(V)=\{\int v(z)z^{m}, v\in V, m\in \BZ\}$.
By (I.1.3.7), $\Lie(V)$ is a Lie subalgebra of $\End(V)$. In particular,
$\Lie(V)$ is a $\Vir$-module and $ad L_{0}$ is diagonalizable:

$$ \Lie(V)=\oplus_{m\in\BZ}\Lie(V)_{m},\
\Lie(V)_{m}=\{A\in\Lie(V):\ [L_{0},A]=mA\}.
\eqno{(1.1.1)}
$$

Observe
that $id\in\End(V)$ belongs to $\Lie(V)$ as $id=\int|0\rangle (z)z^{-1}$.
Below we shall sometimes quotient the Lie ideal $\BC id$ out.

Define:
$$
\Lie(V)_{\leq}=\{A\in\Lie(V):\ (ad L_{-1})^{m}A=0 \text{ if } m>>0\}/\BC id,
\eqno{(1.1.2)}
$$

$$
\Lie(V)_{\leq n}=\oplus_{m\leq n}\Lie(V)_{m}
\eqno{(1.1.3)}
$$

$$
\tlie(V)=\lim_{-\infty\leftarrow m} \Lie(V)/\Lie(V)_{\leq m}
\eqno{(1.1.4)}
$$

$$
\tlie(V)_{\leq}=(\lim_{-\infty \leftarrow m}\Lie(V)_{\leq}/
(\Lie(V)_{\leq m}\cap\Lie(V)_{\leq}))/\BC id.
\eqno{(1.1.5)}
$$

One easily checks that the bracket on $\Lie(V)$ extends by continuity
to that on $\tlie(V)$ making the latter into a Lie algebra,  and that
$\tlie(V)_{\leq}\subset\tlie(V)$ is a Lie subalgebra.

Further, the Lie algebra of vector fields on the formal disk
$\BC[[z]]d/dz\subset\Vir$ operates on $\tlie(V)_{\leq}$, the action
of the subalgebra $z\BC[[z]]d/dz$ being integrable.
       The formal geometry of Gelfand and Kazhdan produces then
a universal sheaf of Lie algebras. It means that for any smooth curve
$\CC$ we get a sheaf of Lie algebras, $\lie_{\CC}(V)$ on $\CC$, so that
$\Gamma(\Spec\BC[[z]], \lie_{\Spec\BC[[z]]}(V))=\tlie(V)_{\leq}$. The
correspondence $\CC\mapsto  \lie_{\CC}(V)$ is functorial with respect
to etale morphisms.

{\bf 1.2.} Here is a somewhat more explicit description of the
sheaf $\lie_{\CC}(V)$. Let $v\in V_{\Delta}$ be a $\Vir$-singular vector.
 This means that $v$ does not belong to
$\BC |0\rangle$ and is annihilated by $L_{i},\ i\geq 1$. In the language
of fields we get, due to (I.1.7.3), that:

$$
L(z)v(w)\sim \frac{v(w)'}{z-w}+\frac{\Delta v(w)}{(z-w)^{2}},
\eqno{(1.2.1)}
$$
that is to say, that $v(z)$ is a {\it primary} field.
A glance  at (I.1.3.7) or (I.1.7.4) shows that
$$
[L_{m}, v_{(n)}]= (-n+(m+1)(\Delta-1))v_{(m+n)}.
\eqno{(1.2.2)}
$$

In other words, under the action of $\Vir$ elements $v_{(n)}\in \Lie(V)$
transform as if $v_{(n)}$ were equal to $z^{n}(dz)^{-\Delta+1}$ and
therefore $\Sigma_{n\geq 0}c_{n}v_{(n)}\in\tlie(V)_{\leq}$ can be regarded
as a formal jet of the section of the bundle of $-\Delta+1$-differentials.
 It
follows that when ``spread'' over a curve a $\Vir$-singular vector
 of conformal weight $\Delta$ gives a subsheaf of $\lie_{\CC}(V)$
isomorphic to the sheaf of $-\Delta+1$ differentials.

Of course not any element of $V$ is singular.
 However $V$ possesses the filtration by conformal weights:
$...V_{\leq 0}\subset V_{\leq 1}\subset..., $
where $V_{\leq m}=\oplus_{i\leq m}V_{i}$.
By definition elements $L_{m}, m>0$ preserve this filtration and
 act trivially on $\oplus_{m}V_{\leq m}/\oplus V_{\leq m-1}$.
Using (I.1.7.4) once again we get that $\tlie(V)_{\leq}$ has a filtration
 such
that the coorresponding graded object is isomorphic, as a
 $\BC[[z]]d/dz$-module,
to a direct sum of modules of $\Delta$-differentials over $\BC^{\ast}$.
There is one such module for each $v\in V_{-\Delta+1}$ if $\Delta\neq 1$
or for each $v\in V_{0}/\BC |0\rangle$. If $v=L_{-1}w,\; v,w\in V$,
then the corresponding modules are equal (as subspaces of $\tlie(V)_{\leq}$).

Putting all of this together we get that
$\lie_{\CC}(V)$ has a filtration
such that the corresponding graded object is a direct sum of sheaves
of $\Delta$-differentials with appropriate $\Delta$.
Further, locally in the presence of a coordinate,
for example when $\CC$ is either $\Spec\BC[[z]]$ or
$\Spec\BC[z,(z-z_{1})^{-1},...,(z-z_{m})^{-1}]$,
$ \lie_{\CC}(V)$ is a
free sheaf of $\CO_{\CC}$-modules whose fiber is isomorhic to
$V/(L_{-1}V+\BC |0\rangle)$.

{\bf 1.3.}  Here we collect several Lie algebras that arise in the case
of a curve with marked points.
By construction $\Lie(V)$ is a central extension of
$\Gamma(\BC^{\ast},\lie_{\BC^{\ast}}(V))$. Likewise
$\tlie(V)$ is a central extension of
$\Gamma(\Spec\BC((z)),\lie_{\Spec\BC((z))}(V))$.

Let $P$ be a point of $\CC$, $U$ the formal neighborhood of $P$,
$U'=U-P$. We have then the Lie algebras $\Gamma(U',\lie_{U'}(V))$
and $\Gamma(\CC-P,\lie_{\CC}(V))$ along with the ``localization at
$P$'' map
$$
i_{P}: \Gamma(\CC-P,\lie_{\CC}(V))\rightarrow\Gamma(U',\lie_{U'}(V)).
\eqno{(1.3.1)}
$$

 Since
any isomorphism $\BC[[z]]\rightarrow \CO_{U}$
(a choice of a local coordinate) determines an isomorphism
$\Gamma(\Spec\BC((z)),\lie_{\Spec\BC((z))}(V))\rightarrow
\Gamma(U',\lie_{U'}(V))$,
 there arises a  Lie algebra,
to be denoted by $\Lie(V)^{P}$, a central extension of
$\Gamma(U',\lie_{U'}(V))$ with a distinguished central element
$\BK$. Again
any isomorphism $\BC[[z]]\rightarrow \CO_{U}$
 determines an isomorphism $\tlie(V)\rightarrow \Lie(V)^{P}$
preserving $\BK$.
The algebra $\Lie(V)^{P}$ is sometimes referred to as
``$\tlie(V)$ sitting at $P$''.
Given a collection of points $\{P_{1},...P_{m}\}\subset\CC$
 consider the Lie algebra $\oplus_{i=1}^{m}\Gamma(U_{i}',\lie_{U_{i}'}(V))$,
where $U_{i}'$ is the formal punctured neighborhood of $P_{i}$,
and its central extension $\Lie(V)^{P_{1},...,P_{m}}$, the Baer sum
of central extensions
$$
0\rightarrow\BC\rightarrow \Lie(V)^{P_{i}}\rightarrow
\Gamma(U_{i}',\lie_{U_{i}'}(V))\rightarrow 0
$$
(To take the  Baer sum
of central extensions by $\BC$ means to take the direct
sum of central extensions
and then set $\BK$ in all summands equal each other.)

Another ingredient is the Lie algebra
$\Gamma(\CC-\{P_{1},...,P_{m}\},\lie_{\CC}(V))$; because of its
importance we shall denote it simply as $\Lie(V)_{out}$, suppressing
some of the data on which it actually depends. There is a Lie
algebra morphism (cf. (1.3.1))
$$
i_{P_{1},...P_{m}}=\oplus_{j} i_{P_{j}}:
 \Lie(V)_{out}\rightarrow
\oplus_{i=1}^{m}\Gamma(U_{i}',\lie_{U_{i}'}(V)).
\eqno{(1.3.2)}
$$

The key point is that the last map lifts to a Lie algebra morphism
$$
\hat{i}_{P_{1},...P_{m}}=\oplus_{j} i_{P_{j}}:
 \Lie(V)_{out}\rightarrow   Lie(V)^{P_{1},...,P_{m}}
\eqno{(1.3.3)}
$$

This can be proved by methods of [BFM, 2.2.1-2.3.4], that is, by
showing that the  central extensions in question are Tate ones and then
using the residue theorem. We shall skip this argument as below
in the case of interest for us we shall exhibit a direct construction.

{\bf 1.4.} Let us now discuss in greater detail the case when
$\CC=\cp$. We assume fixed an embedding $\BC\hookrightarrow\cp$ and
a coordinate $z$ on $\BC$. First of all we shall explain that
to construct
the splitting (1.3.3) in the case of $\cp$ means to construct
  a certain Lie algebra antiinvolution
$\eta: \Lie(V)\rightarrow\Lie(V)$ preserving the central element.

Begin with the case of 2 marked points  and let the points be
$0$ and $\infty$. In this case $\Lie(V)_{out}=
\Gamma(\BC^{\ast},\lie_{\BC^{\ast}}(V))$. Observe that the coordinate
change $z\mapsto 1/z$ induces a Lie algebra involution:
$$
\bar{\eta}:
\Gamma(\BC^{\ast},\lie_{\BC^{\ast}}(V))\rightarrow
\Gamma(\BC^{\ast},\lie_{\BC^{\ast}}(V)).
\eqno{(1.4.1)}
$$

 Choose the local coordinates at $0$ and $\infty$
to be $z$ and $1/z$ respectively. This allows to identify
$\Gamma(U_{0}',\lie_{U_{0}'}(V))=
\Gamma(U_{\infty}',\lie_{U_{\infty}'}(V))=\tlie(V)/\BC\BK$. Under
this identification, both $i_{0}$ and $i_{\infty}$ are
 isomorphisms of $\Lie(V)_{out}$
with $\Lie(V)/\BC\BK\subset\tlie(V)/\BC\BK$ such that
$i_{\infty}=i_{0}\circ\bar{\eta}$. We can and will assume that
$i_{0}= id$ (it is simply a matter of notation); then
$i_{\infty}=\bar{\eta}$.

By definition both $i_{0}$ and $i_{\infty}$ lift to isomorphisms
of central extensions
$$
\hat{i}_{0},\hat{i}_{\infty}:\Lie(V)\rightarrow\Lie(V).
\eqno{(1.4.2)}
$$
Because of our conventions, we can set $\hat{i}_{0}=id$.
Then existence of the splitting (1.3.3) in the situation
in question is equivalent
to existence of $\hat{i}_{\infty}$ in (1.4.2) so that
$
\hat{i}_{\infty}(\BK)=-\BK.
$
 It follows that the map $\eta$
defined to be equal $-\hat{i}_{\infty}$ satisfies the following
conditions:
$$
\eta:\Lie(V)\rightarrow\Lie(V) \text{ is an antiinvolution},
\eqno{(1.4.3)}
$$
$$
\eta(\BK)=\BK,
\eqno{(1.4.4)}
$$
$$
\eta |_{\Lie(V)/\BC\BK}=-\bar{\eta}.
\eqno{(1.4.5)}
$$

As we explained, the data (1.4.3-5) is equivalent to the datum
(1.3.3) in the case of 2 marked points $0,\infty\in\cp$.

{\bf 1.5.} Let us now construct the spliting (1.3.3) in the case of several
points on $\cp$. For the sake of definiteness assume that $P_{1}=\infty$,
while $P_{i}=z_{i}\in\BC,\; i=2,,,,,m$. Choose $z-z_{i}$ as a local
coordinate at $P_{i}$ and $1/z$ as that at $\infty$. Trivializing the
sheaf $\lie_{\cp}(V)$ over $\BC$ using the coordinate $z$ (see 1.2) we
identify
the Fourier
component $B_{(j)}$ of any field $B(z)$ with a section
of a trivial bundle over $\BC$ with fiber: image of
 $\BC B_{(-1)}|0\rangle\subset V$
in $V/(L_{-1}V+\BC|0\rangle)$, see the very end of 1.2.
  The translation invariance axiom
shows then that under this identification $B_{(j)}$
becomes $z^{j}$. Similarly,
denote by $B_{(j)}^{z_s}$ the section of the same bundle
 identified with the function $(z-z_{s})^{j}$; in particular,
$B_{j}^{0}=B_{j}$. In this
notation the localization maps (1.3.1-2) are simply Laurent series
expansions
of rational functions. We have:
$$
i_{P_{1},...,P_{m}}:
\Lie(V)_{out}\rightarrow\oplus_{t=1}^{m}\Gamma(U_{t}',\lie_{U_{t}'}(V)),
 $$
$$
i_{P_{1},...,P_{m}}=\oplus_{t}i_{P_{t}},
\eqno{(1.5.1)}
$$
where

$$
i_{P_{t}}: \Lie(V)_{out}\rightarrow
\Gamma(U_{t}',\lie_{U_{t}'}(V)),
\eqno{(1.5.2)}
$$
is defined as follows

$$
i_{z_{t}}(B_{(j)}^{z_{s}})=
-\partial_{z_{s}}^{(-j+1)}(\Sigma_{i=0}^{\infty}
\frac{B_{(i)}^{z_{t}}}{(z_{s}-z_{t})^{i+1}}) \text{ if } j<0,
\eqno{(1.5.3)}
$$
(note that the expression in the r.h.s. of the last equality is
nothing but $(B_{(j)}|0\rangle)(z_{s}-z_{t})_{-}$, see (I.1.4.3) for
the definition)

$$
i_{\infty}(B_{(j)}^{z_{s}})=-\eta((B_{(j)}|0\rangle)(z_{s})_{+})
\text{ if }j<0,
\eqno{(1.5.4)}
$$
(see (I.1.4.3) for the definition of $(B_{(j)}|0\rangle)(z_{s})_{+}$)

$$
i_{z_{t}}(B_{(j)}^{z_{s}})=\Sigma_{i\geq 0}\binom{j}{i}(z_{t}-z_{s})^{j-i}
B_{i}^{z_{t}} \text{ if } j\geq 0,
\eqno{(1.5.5)}
$$

$$
i_{\infty}(B_{(j)})=-\eta(B_{(j)} )
\text{ if }j\in\BZ.
\eqno{(1.5.6)}
$$

{\bf Remark.} The formula (1.5.3) is obtained by expanding the function
$1/z-z_{s}$ at $z_{t}$:
$$
\frac{1}{z-z_{s}}=-\Sigma_{i=0}^{\infty}\frac{(z-z_{t})^{i}}
{(z_{s}-z_{t})^{i+1}}
$$
and identifying $(z-z_{t})^{i}$ in the r.h.s. with $B^{z_{t}}_{(i)}$.
A similar remark applies to (1.5.5) with $1/z-z_{s}$ replaced
with $(z-z_{s})^{j}$. Of course (1.5.3) follows from (1.5.5).

\bigskip

In the presence of a coordinate at any point $P$
 any splitting $V=V'\oplus\BC |0\rangle$ gives a splitting $\Lie(V)^{P}=
\Gamma(U,\lie_{U}(V))\oplus\BC\BK,\; P\in U$.
We have already fixed coordinates at our points and we now choose
arbitrarily a splitting $V=V'\oplus\BC |0\rangle$ to get splittings
$$
\Lie(V)^{P_{t}}=\Gamma(U_{t}',\lie_{U_{t}'}(V))\oplus\BC\BK.
\eqno{(1.5.7)}
$$
 Given that, there
arise canonical embeddings
$\Gamma(U_{t}',\lie_{U_{t}'}(V))\hookrightarrow \Lie(V)^{P_{t}}$.
We use these embeddings to define (cf. (1.3.2-3)):
$$
\hat{i}_{P_{1},...,P_{m}}:\Lie(V)_{out}\rightarrow\oplus_{t=1}^{m}
\Lie(V)^{P_{t}}
$$
to be
$$
\hat{i}_{P_{1},...,P_{m}}=i_{P_{1},...,P_{m}}.
\eqno{(1.5.8)}
$$

In the case of 2 points, $\infty$ and $0$ the fact that this
$\hat{i}_{P_{1},...,P_{m}}$ is a Lie algebra homomorphism follows
from 1.4. If the points in question are $\infty$ and $z_{2}$  then
the desired
claim follows from the following observation: operator $L_{-1}$
annihilates the cocycle determined by the splittings
(1.5.7). Finally, the case of more than 2 points is readily reduced
 to that of 2 points by representing a rational function as
 a linear combination of $(z-z_{t})^{j}$.

{\bf 1.6.}  So far we have been dealing with ``algebras attached to
points on a curve''. One can as well attach to points modules over these
algebras. More precisely, let $\CC$ be a curve with $m$ marked points
$P_{1},...,P_{m}$ and coordinates $z_{i}$ on the formal neighborhoods
$U_{i}\ni P_{i},\; i=1,...,m$. This gives  canonical identifications
$\Lie(V)=\Lie(V)^{P_{i}},\; i=1,...,m$. Therefore, if $M_{1},...,M_{m}$
is a collection of $\Lie(V)$-modules of the same central charge, then
$M_{1}\otimes M_{2}\otimes\cdots\otimes M_{m}$ is a $\Lie(V)^{P_{1},...,
P_{m}}$-module, meaning that elements of $\Lie(V)^{P_{i}}$ act on
the $i$-th factor of the tensor product
$M_{1}\otimes M_{2}\otimes\cdots\otimes M_{m}$.

Due to (1.3.3),  $M_{1}\otimes M_{2}\otimes\cdots\otimes M_{m}$
is also a $\Lie(V)_{out}$-module. Hence there arises the space
of coinvariants
$$
(M_{1}\otimes M_{2}\otimes\cdots\otimes M_{m})_{\Lie(V)_{out}}:=
\frac{M_{1}\otimes M_{2}\otimes\cdots\otimes M_{m}}
{\Lie(V)_{out}(M_{1}\otimes M_{2}\otimes\cdots\otimes M_{m})}.
\eqno{(1.6.1)}
$$

We shall often make use of the dual space
$$
\langle M_{1},M_{2},..., M_{m}\rangle(P_{1},...,P_{m}) :=
 ((M_{1}\otimes M_{2}\otimes\cdots\otimes M_{m})_{\Lie(V)_{out}})^{\ast}.
\eqno{(1.6.2)}
$$

Explicitly, $\langle M_{1},M_{2},..., M_{m}\rangle(P_{1},...,P_{m})$
 consists of linear functionals
$$
\langle.,.,...,.\rangle: M_{1}\otimes M_{2}\otimes\cdots\otimes M_{m}
\rightarrow\BC,\; v_{1}\otimes\cdots v_{m}\mapsto
\langle v_{1},...,v_{m}\rangle
$$
satisfying the following $\Lie(V)_{out}$-invariance condition:
$$
\Sigma_{t=1}^{m}\langle  v_{1},...,v_{t-1},i_{P_{t}}(X)v_{t},v_{t+1},...
v_{m}\rangle=0
\eqno{(1.6.3)}
$$
for any $X\in\Lie(V)_{out}$.

By definition the space $\langle M_{1},M_{2},..., M_{m}\rangle
(P_{1},...,P_{m})  $
depends on the choice of local coordinates. However, if the modules
$M_{i}$ are all integrable with respect to the algebra $\BC[[z]]d/dz$,
then the two such spaces associated with different coordinates are
canonically identified. This means that there arises a vector bundle
$\langle M_{1},M_{2},..., M_{m}\rangle_{\CC}$
over the configuration space
$\CC^{\times m}-\text{ diagonals }$,
the fiber over the point $(P_{1},...,P_{m})$ being
 $\langle M_{1},M_{2},..., M_{m}\rangle(P_{1},...,P_{m})$.

Denote by $L_{-1}^{(i)}$ the linear transformation of
$\langle M_{1},M_{2},..., M_{m}\rangle(P_{1},...,P_{m})$ acting
 as $L_{-1}$ on the
$i$-th factor.
It is well known
that the operators
$$
\Delta_{t}=\partial/\partial z_{t} -L_{-1}^{t},\; t=1,...,m
\eqno{(1.6.4)}
$$
define a flat connection on the bundle
$\langle M_{1},M_{2},..., M_{m}\rangle_{\CC}$.

We now try to explain that
in many cases this construction admits a transparent
 representation-theoretic or ``vertex-theoretic'' interpretation.

{\bf 1.7.} Denote by $V-\Mod$ the category of restricted graded $V$-modules.
Existence of the antiinvolution $\eta$ (see (1.4.3)) allows us to define
the following duality functor
$$D: V-\Mod\rightarrow V-\Mod,\; M\mapsto D(M):=
\oplus_{n\in\BZ}M_{n}^{\ast}.
\eqno{(1.7.1)}
$$
 The $V$-module structure on $D(M)$ is defined in the following manner:
$$
(b(z)f)(.):=f(\eta(b(z)).)=\Sigma_{i\in\BZ}z^{-i-1}f(\eta(b_{(i)}.))
\eqno{(1.7.2)}
$$
$$\text{ for all } b\in V, f(.)\in D(M).
$$

Now consider the case of the 2 marked points,
$\infty$ and $0$, on $\cp$ with the standard choice of local coordinates --
$z$ and $ z^{-1}$. Comparing the definition of $\eta$ in 1.3, or
the formulas in 1.5, with (1.6.3)  we immediately get the canonical
embedding:
$$
Hom_{V}(M_{1},M_{2})\hookrightarrow
\langle D(M_{2}), M_{1}\rangle(\infty,0),.
\eqno{(1.7.3)}
$$
$$
F\mapsto\{(x,y)\mapsto y(F(x))\}.
$$

Further, if homogeneous subspaces $M_{2}$ are finite dimensional, then
(1.7.3) is an isomorphism. This finiteness condition sometimes fails as
it does in some of our concrete examples to be considered below.

{\bf 1.8.} Consider now the case of 3 points  on $\cp$.
Attach modules $D(V), V, V$ to the points $\infty$, $z$ and
$0$  respectively.

{\bf Lemma.} {\it The functional $\Phi_{z}(.,.,.)\in
(D(V)\otimes V\otimes V)^{\ast}$ defined by the formula
$$
\Phi_{z}(a,b,c)=a(b(z)c)
\eqno{(1.8.1)}
$$
actually belongs to $\langle D(V),V,V\rangle(\infty,z,0)$.}

\bigskip

 Consider $\langle D(V),V,V\rangle_{\BC^{\ast}}$, a vector
bundle over $\BC^{\ast}$ obtained by restricting
$\langle D(V),V,V\rangle_{\CC}$, a trivial bundle over
$(\cp)^{\times 3}-\text{ diagonals }$ defined in 1.6, to
$\{(\infty,z, 0),\; z\in\BC^{\ast}\}\subset
(\cp)^{\times 3}-\text{ diagonals}$. Of course, the connection
(1.6.4) restricts to the $\langle D(V),V,V\rangle_{\BC^{\ast}}$.
 The translation invariance axiom (see I.1.2) along
with (I.1.7.4) implies that the section
$z\mapsto\Phi_{z}$ is horizontal. Indeed we have:
$$
\Delta (\Phi_{z})(a,b,c)= a((b(z)'-(L_{-1}b)(z))c)
$$
$$
= a((b(z)'-[L_{-1},b(z)])c)=0.
$$

{\bf 1.9.} More generally, we have the following result. Observe
that although the formal product of fields
 $$
a_{1}(z_{1})a_{2}(z_{2})a_{3}(z_{3})\cdots
a_{m-1}(z_{m-1})a_{m}(z_{m}),\; a_{i}\in V
$$
does not make sense as an operator on $V$, the matrix element
$$
y (a_{1}(z_{1})a_{2}(z_{2})a_{3}(z_{3})\cdots
a_{m-1}(z_{m-1})a_{m}(z_{m})x),\; x\in V,y\in D(V)
$$
is a well-defined formal
Laurent series. Further, it easily follows from the OPE formula
(I.1.6.9) that this series converges to a certain rational function
in the region $|z_{1}|>|z_{2}|>\cdots >|z_{m}|$.

{\bf Lemma.} {\it The functional $\Phi_{z_{2},...,z_{m-1}}(.,...,.)\in
(D(V)\otimes V\otimes\cdots\otimes V)^{\ast}$ defined by the formula
$$
\Phi_{z_{2},...,z_{m-1}}(a_{1},a_{2},...,a_{m-1},a_{m})=
a_{1}(a_{2}(z_{2})a_{3}(z_{3})\cdots a_{m-1}(z_{m-1})a_{m})
\eqno{(1.9.1)}
$$
actually belongs to $\langle D(V),V,...,V\rangle(\infty,z_{2},...
z_{m-1},0)$.}

{\bf Proof.} By way of preparation let us remark that due to the locality
axiom the OPE (I.1.4.12) formula
 can be rewritten as the following 2 identities
(see also [K], Theorem 2.3 (ii)):
$$
[a_{s}(z_{s}),B(z_{t})_{+}]=\Sigma_{i=0}^{\infty}
\frac{(B_{(i)}a_{s})(z_{s})}{(z_{t}-z_{s})^{i+1}},
\eqno{(1.9.1)}
$$
$$
[B(z_{t})_{-},a_{s}(z_{s})]=\Sigma_{i=0}^{\infty}
\frac{(B_{(i)}a_{s})(z_{s})}{(z_{t}-z_{s})^{i+1}}.
\eqno{(1.9.2)}
$$

When compared to (1.5.3) the last equalities rewrite as follows:
$$
[a_{s}(z_{s}),B(z_{t})_{+}]= -(i_{z_{s}}(B_{(-1)}^{z_{t}}))(a_{s}))(z_{s}),
\eqno{(1.9.3)}
$$
$$
[B(z_{t})_{-},a_{s}(z_{s})]= -(i_{z_{s}}(B_{(-1)}^{z_{t}}))(a_{s}))(z_{s}).
\eqno{(1.9.4)}
$$

Now turn to the proof proper.
We have to show that
$$
(X\Phi_{z_{2},...,z_{m-1}})(a_{1},a_{2},...,a_{m-1},a_{m})=0
$$
for any $X\in\Lie(V)_{out}$. When restricted to $\BC$, $X$ may only
have poles at $z_{t},\; 2\leq t\leq m-1$. It means that it is enough
to consider the following cases: $X$ is either $B^{z_{t}}_{(j)},\; j<0$
or $B_{(j)},\; j\geq 0$. Let for simplicity $X$ be equal
$B_{(-1)}^{z_{t}},\; 1<t<m$. We have:
$$
a_{1}(a_{2}(z_{2})\cdots a_{t-1}(z_{t-1})(B_{(-1)}a_{t})(z_{t})
a_{t+1}(z_{t+1}),...
,a_{m-1}(z_{m-1}a_{m})
$$
$$
=
a_{1}(a_{2}(z_{2}),...,B(z_{t})_{+}a_{t}(z_{t}),...
,a_{m-1}(z_{m-1})a_{m})
$$
$$
+ a_{1}(a_{2}(z_{2}),...,a_{t}(z_{t})B(z_{t})_{-},...
,a_{m-1},a_{m})
$$
$$
=\Sigma_{s=2}^{t-1}
 a_{1}(a_{2}(z_{2}),...,[a_{s}(z_{s}),B(z_{t})_{+}],...,a_{t}(z_{t}),...
,a_{m-1}(z_{m-1})a_{m})
$$
$$
+
\Sigma_{s=t+1}^{m-1} a_{1}(a_{2}(z_{2}),...,a_{t}(z_{t}),...,
[B(z_{t})_{+},a_{s}(z_{s})],...,
,a_{m-1}(z_{m-1})a_{m})
$$
$$
+a_{1}(B(z_{t})_{+}a_{2}(z_{2}),...,a_{t}(z_{t}),...
,a_{m-1}(z_{m-1})a_{m})
$$
$$
+
a_{1}(a_{2}(z_{2}),...,a_{t}(z_{t}),...
,a_{m-1}(z_{m-1})B(z_{t})_{-}a_{m}).
$$

The terms in the r.h.s of this equality are interpreted as follows:

 by (1.9.3,4)
 each summand in the summations $\Sigma_{s}$ equals:
$$
-a_{1}(a_{2}(z_{2}),...,(i_{z_{s}}(B^{z_{t}}_{(-1)}a_{s})(z_{s}),...
,a_{m-1}(z_{m-1})a_{m});
$$

the last summand equals, directly by (1.5.3),
$$
-a_{1}(a_{2}(z_{2}),...,
,a_{m-1}(z_{m-1})(i_{z_{m}=0}(B^{z_{t}}_{(-1)})a_{m})),
$$

while the one before the last one, by (1.5.4), equals
$$
-i_{z_{1}=\infty}(B^{z_{t}}_{(-1)})a_{1}(a_{2}(z_{2})...
a_{m-1}(z_{m-1})a_{m}).
$$

Collecting all the terms in the l.h.s. we obtain the desired equality
$B^{z_{t}}_{(-1)}\Phi_{z_{2},...,z_{m-1}}=0$. The case of $B^{z_{t}}_{j},
j<-1$ is obtained from the one considered by taking derivatives. The
case $X=B_{(j)}, j\geq 0$ is treated similarly.

\bigskip

As in 1.8 one can consider the embedding
$$
\pi:(\BC^{\ast})^{\times(m-2)}-\{\text{ diagonals }\}\hookrightarrow
\CC^{\times m}-\{\text{ diagonals }\}
\eqno{(1.9.5)}
$$
obtained by keeping two points equal 0 and $\infty$. There arises the
pull-back
$\pi^{\ast}\langle D(V),V,...,V\rangle$ to $(\BC^{\ast})^{\times(m-2)}$
and one proves, exactly as in the end of 1.8,
that $(z_{2},...,z_{m-1})\mapsto \Phi_{z_{2},...,z_{m-1}}(.,...,.)$
is its horizontal section.

\bigskip

\bigskip

\centerline{\bf   \S 2.   An application to the chiral De Rham complex}

\bigskip

Here we shall explain what the abstract constructions of \S 1 mean
in several concrete examples. In the end we shall explain how to recover
the multiplication in the cohomology ring of  a smooth manifold
by looking at correlation functions.

{\bf 2.1.} Consider the vertex algebra $\CO_{N}$. It is generated by
the fields

$$
a^{i}(z)=\Sigma_{s\in\BZ}a^{i}_{s}z^{-s-1}
\eqno{(2.1.1)}
$$
$$
b^{i}(z)=\Sigma_{s\in\BZ}b^{i}_{s}z^{-s}
\eqno{(2.1.2)}
$$
so that the Fourier components $a^{i}_{s}, b^{j}_{t}\in\Lie(\CO_{N})$
satisfy the relations
$$
[a^{i}_{s},b^{j}_{t}]=\delta_{ij}\delta_{s,-t}\BK.
\eqno{(2.1.3)}
$$
These relations mean that we are dealing with
 the Heisenberg algebra, $H_{N}$, linearly
spanned by $a$'s and $b$'s. $\CO_{N}$ is a level 1 ($\BK\mapsto 1$)
 representation
of this algebra; it is generated by the ``vacuum'' vector, $|0\rangle$,
satisfying the relations:
$$
a^{i}_{s}|0\rangle=b^{i}_{t}|0\rangle=0\text{ if } s\geq 0,t>0.
$$
  The Virasoro
field is as follows:
$$
L(z)=\Sigma_{i=1}^{N}:b^{i}(z)'a^{i}(z):
\eqno{(2.1.4)}
$$
One easily checks that $b^{i}_{0}|0\rangle$ is a $\Vir$-singular vector
of weight 0. It follows that (cf. 1.2):
$$
[L_{s},a^{i}_{j}]=-j a^{i}_{s+j},
\eqno{(2.1.5)}
$$
$$
[L_{s},b^{i}_{j}]=-(s+j)b^{i}_{s+j}.
\eqno{(2.1.6)}
$$
It follows that as elements of $\Gamma(\BC^{\ast},\lie_{\BC^{\ast}}
(\CO_{N}))$, $a^{i}_{j}$ and $b^{i}_{j}$ are identified with
$z^{j}$ and $z^{j-1}dz$ respectively. Hence the automorphism
$\bar{\eta}$ of $\Gamma(\BC^{\ast},\lie_{\BC^{\ast}}
(\CO_{N}))$
 induced by the coordinate change $z\mapsto z^{-1}$
(see (1.4.1)) operates on these elements as follows:
$$
\bar{\eta}(a^{i}_{j})=a^{i}_{-j},\;
\bar{\eta}(b^{i}_{j})=-b^{i}_{-j}.
\eqno{(2.1.7)}
$$
These formulas suggest how to lift $\bar{\eta}$  to an automorphism
$\hat{i}_{\infty}$ (cf. (1.4.2) ) at least when restricted to $H_{N}$:
$$
\hat{i}_{\infty} (a^{i}_{j})=a^{i}_{-j},\;
\hat{i}_{\infty}(b^{i}_{j})=-b^{i}_{-j},\
\hat{i}_{\infty}(\BK)=-\BK
\eqno{(2.1.8)}
$$

It is now easy to find the antiautomorphism $\eta$; here  is how it operates
on $H_{N}$:
$$
\eta(a^{i}_{j})=-a^{i}_{-j},\;
\eta(b^{i}_{j})=b^{i}_{-j},\; \eta(\BK)=\BK.
\eqno{(2.1.9)}
$$

It is also easy to extend the {\it anti}automorphism
 $\eta$ from $H_{N}$ to the entire $\Lie(\CO_{N})$ -- easier
 than the automorphism $\hat{i}_{\infty}$.
 Each element of $\Lie(V)$ is an infinite
sum of monomials in $a$'s and $b$'s; to evaluate $\eta$ on such a series
one has to apply $\eta$ to each summand regarded as an element of the
universal enveloping algebra of $H_{N}$. A similar procedure does not apply
to $\hat{i}_{\infty}$.

{\bf 2.2.} As we argued in [MSV], the vertex algebra $\CO_{N}$ is associated
with $\BC^{N}$ with a fixed coordinate system; this is
reflected, in particular, in the fact that the conformal weight
0 component of $\CO_{N}$ is $\BC[b^{1}_{0},...,b^{N}_{0}]$.
 Passing to various completions we get the vertex algebras $\hat{\CO}_{N}$,
$\CO_{N}^{an}(U)$, $\CO_{N}^{sm}(U)$ associated with a formal disk,
  an open set
in analytic category and an open set in smooth category resp.; as the
conformal weight 0 component they have respectively: the algebra of formal
power series in $b^{1}_{0},...,b^{N}_{0}]$,
the algebra of analytic functions, or the algebra of
smooth functions on the given open set $U$.
What was said in 2.1 carries over
to these completions without serious changes.

Let $f(b)$ be a non-constant
 function in any of the mentioned categories identified
with an element of the conformal weight 0 component of the corresponding
vertex algebra $\CO_{N}^{\cdot}$. Denote also by $f(b)$ the
corresponding element
 $\CO^{\cdot}_{N}$ and by $f(b)_{j},\; j\in\BZ$
the corresponding elements of
 $\Gamma(\BC^{\ast},\lie_{\BC^{\ast}}(\CO_{N}^{\cdot})$.

The Virasoro field is given by (2.1.4). As $f(b)|0\rangle$ is a
$\Vir$-singular vector of weight 0 we get
(analogously to (2.1.6):

$$
[L_{s},f(b)_{j}]=-(s+j)f(b)_{s+j}.
\eqno{(2.2.1)}
$$

Therefore we have the following analogue of (2.1.9)
$$
\eta(a^{i}_{j})=-a^{i}_{-j},\;
\eta(f(b)_{j})=f(b)_{-j},\; \eta(\BK)=\BK.
\eqno{(2.2.2)}
$$

An extension of (2.2.2) to the entire $\Lie(\CO_{N}^{an})$ or
$\Lie(\CO_{N}^{sm})$ is exactly as explained in the end of 2.1.

Observe by the way that in any restricted vertex algebra $V$
Fourier components of a field  $v(z)$ associated with $v\in V_{0}$
  transform as differential forms.

{\bf 2.3.} Our exposition in \S 1 was suited for the pure even case.
It is obvious however that the same could have been done in the
supercase by changing signs in certain places in the standard way.
Consider, for example, the vertex algebra $\Lambda_{N}$ based
on the Clifford algebra $Cl_{N}$ in the same way as $\CO_{N}$ is
based on $H_{N}$. The basis of $Cl_{N}$ is as follows:
$\phi^{i}_{j},\psi^{i}_{j}$ (odd), $\BK$ (even); the relations are:
$$
[\phi^{i}_{s},\psi^{j}_{t}]=\delta_{ij}\delta_{s,-t}\BK
\eqno{(2.3.1)}.
$$
The vacuum $|0\rangle\in\Lambda_{N}$ satisfies
$$
\phi^{i}_{s}|0\rangle=\psi^{j}_{t}|0\rangle=0,\; s>0,t\geq 0,
\eqno{(2.3.2)}
$$
$$
\BK|0\rangle=|0\rangle.
\eqno{(2.3.3)}
$$

The Virasoro field is given by
$$
L(z)=\Sigma_{i=1}^{N}:\phi^{i}(z)'\psi^{i}(z):
\eqno{(2.3.4)}
$$
Similarly to (2.1.5-6) we get
$$
[L_{s},\psi^{i}_{j}]=-j \psi^{i}_{s+j},
\eqno{(2.3.5)}
$$
$$
[L_{s},\phi^{i}_{j}]=-(s+j)\phi^{i}_{s+j}.
\eqno{(2.3.6)}
$$

It follows that $\phi^{i}_{j}$ is identified with $z^{j-1}dz$, while
$\psi^{i}_{j}$ -- with $z^{j}$. Going through the same steps is in 2.1,
we obtain the following formulas for the antiinvolution $\eta$ (cf. (2.1.9):

$$
\eta(\psi^{i}_{j})=-\psi^{i}_{-j},\;
\eta(\phi^{i}_{j})=\phi^{i}_{-j},\; \eta(\BK)=\BK.
\eqno{(2.3.7)}
$$

(To check that the relations (2.3.1) are indeed preserved one has
to really use the sign rule.) The recipe to extend $\eta$ to the entire
$\Lie(\CO_{N})$, which is explained in the sentences following (2.1.9),
carries over to the present situation word for word.

{\bf 2.4.} What was done in 2.1-2.2 and 2.3. can be combined
and applied in the obvious manner to the vertex algebras
$\Omega_{N}:=\CO_{N}\otimes\Lambda_{N}$,
$\hat{\Omega}_{N}:=\hat{\CO}_{N}\otimes\Lambda_{N}$,
$\Omega_{N}^{an}:=\CO_{N}^{an}\otimes\Lambda_{N}$,
$\Omega_{N}^{sm}:=\CO_{N}^{sm}\otimes\Lambda_{N}$.  An $\eta$ is defined
to act on $a$'s and $b$'s by (2.1.9) and on $\phi$'s and $\psi$'s by
(2.3.7).

All the algebras we just listed possess the four remarkable fields
$$
L(z)=\Sigma_{i=1}^{N}(:b^{i}(z)'a^{i}(z):+:\phi^{i}(z)'\psi^{i}(z):),
\eqno{(2.4.1)}
$$
(the Virasoro field),
$$
G(z)=\Sigma_{i=1}^{N}:b^{i}(z)'\psi^{i}(z):,\;
J(z)=\Sigma_{i=1}^{N}:\phi^{i}(z)\psi^{i}(z):,
Q(z)=\Sigma_{i=1}^{N}:\phi^{i}(z)a^{i}(z):
$$
$$
\eqno{(2.4.2)}
$$

Fourier components of these fields satisfy the
 $N=2$-superconformal algebra relations. We shall not list all of these
here restricting ourselves to the following ones:
$$
[L_{i},Q_{j}]=-jQ_{i+j},
\eqno{(2.4.3a)}
$$
$$
[Q_{0},G(z)]=-L(z).
\eqno{(2.4.3b)}
$$
According to (1.2), (2.4.3a) means that $Q_{-1}|0\rangle$ generates a
subsheaf of $\lie_{\CC}(\Omega_{N}^{\cdot})$
 isomorphic to the structure sheaf for any smooth curve $\CC$.
In particular,
$$
Q_{0}=\int Q(z)\in\Gamma(\CC,\lie_{\CC}(\Omega_{N}^{\cdot})),
\eqno{(2.4.4)}
$$
since $Q_{0}$ is identified with a constant function.

{\bf 2.5.} We now extend the antiinvolution $\eta$ on $\Omega_{N}^{\cdot}$
(see (2.1.9, 2.3.7))
to the sheaf of vertex algebras
$\Omega^{ch}_{X}$, which was associated with a smooth manifold $X$
in [MSV].  This can be easily and uniformly
 explained for either an analytic or
$C^{\infty }$-manifold $X$.

Recall that
$\Omega^{ch}_{X}$ was defined  in each of the 2 settings by considering
an atlas of $X$, associating with each chart a copy of a vertex algebra
of the type considered in 2.2
 and then gluing the algebras
over intersections by lifting the usual transition functions to
operators acting on the algebras. Let us write down the relevant  formulas.

Let $U$ and $\tilde{U}$ be two open subsets of $X$, $b^{1},...,b^{N}$ and
$\tilde{b}^{1},...,\tilde{b}^{N}$ being the respective coordinate functions.
On $U\cap\tilde{U}$ they are related by
$$
\tilde{b}^{i}=g^{i}(b^{1},...,b^{N});\;
b^{i}=f^{i}(\tilde{b}^{1},...\tilde{b}^{N}).
\eqno{(2.5.1)}
$$

Associated to $U$ and $\tilde{U}$ there are $\CO_{N}(U)$ and
$\CO_{N}(\tilde{U})$; the first one is generated by the fields
$a^{i}(z),b^{i}(z),\phi^{i}(z),\psi^{i}(z)$, the second one
by the fields
$\tilde{a}^{i}(z),\tilde{b}^{i}(z),\tilde{\phi}^{i}(z),\tilde{\psi}^{i}(z)$.
 We also have 2 vertex algebras associated to
$U\cap\tilde{U}$: one is an extension of $\CO_{N}(U)$ another is that of
$\CO_{N}(\tilde{U})$. We identify these two algebras by assuming that
the fields are related by the following transformation $G$ lifting (2.5.1):
$$
G(\tilde{b}^{i}(z))=g^{i}(b^{1},...,b^{N})(z),
\eqno{(2.5.2)}
$$
$$
G(\tilde{\phi}^{i}(z))=(\frac{\partial g^{i}}{\partial b^{j}}\phi^{j})(z),
\eqno{(2.5.3)}
$$

$$
G(\tilde{\psi}^{i}(z))= (\frac{\partial f^{i}}{\partial \tilde{b}^{j}}
\psi^{j})(z),
\eqno{(2.5.4)}
$$
$$
G(\tilde{a}^{i}(z))=a^{j}(\frac{\partial f^{i}}{\partial \tilde{b}^{j}})(z) +
(\frac{\partial^{2}f^{k}}{\partial\tilde{b}^{i}\partial\tilde{b}^{l}}
\frac{\partial g^{l}}{\partial b^{r}}\phi^{r}\psi^{k})(z).
\eqno{(2.5.5)}
$$

Denote by $\eta_{U}$ ($\eta_{\tilde{U}}$ resp.)  the antiautomorphism
$\eta$ specialized to $\CO_{N}(U)$
($\CO_{N}(\tilde{U})$ resp.).
It is an easy exercise on definitions in [MSV] to show
that
$$
\eta_{U}\circ G=G\circ\eta_{\tilde{U}}.
\eqno{(2.5.6)}
$$

This equality of course means that the collection of antiinvolutions
$\{\eta_{U}\}$ glues to give  an antiinvolution of the sheaf:
$$
\eta:\ Lie(\Omega^{ch}_{X})\rightarrow Lie(\Omega^{ch}_{X})
\eqno{(2.5.7)}
$$
{\bf 2.6.}
Of the four fields in (2.4.1-2) only two, $G(z)$ and $L(z)$, are
in general invariant
under the transformation (2.5.2-5); it follows, in particular, that
$\Gamma(X,\Omega_{X}^{ch})$ is a restricted conformal vertex algebra.
The two other fields, $Q(z)$ and $J(z)$, are only invariant up to the
addition of  a total derivative; it follows that the Fourier modes
$J_{0}$ and $Q_{0}$ are operators acting on the sheaf $\Omega_{X}^{ch}$.
The action of $J_{0}$ on $\Gamma(X,\Omega_{X}^{ch})$ arising in this
way equips $\Gamma(X,\Omega_{X}^{ch})$ with a gradation (by
eigenspaces of $J_{0}$, or, as they say, by fermionic
charge), while $Q_{0}$ satisfies $Q_{0}^{2}=0$. It follows that
$\Gamma(X,\Omega_{X}^{ch})$ is a complex (infinite in both directions)
 and we proved in [MSV] that its cohomology
$H_{Q}(\Gamma(X,\Omega_{X}^{ch}))$ is canonically isomorphic to the
de Rham cohomology of $X$: $H_{DR}(X)$. The de Rham complex
$C_{DR}(X)$ is, in fact,
canonically identified with the conformal weight 0 component of
$\Gamma(X,\Omega_{X}^{ch})$ and this embedding is a quasiisomorphism.

Given any  $\Gamma(X,\Omega_{X}^{ch})$-module $M$ one can similarly
consider the $Q$-cohomology group
$$
H_{Q}(M)=\frac{Ker (Q_{0}:M\rightarrow M)}{Im (Q_{0}:M\rightarrow M)}.
$$
For example, when applied to $D(\Gamma(X,\Omega_{X}^{ch}))$, this gives
$$
H_{Q}(D(\Gamma(X,\Omega_{X}^{ch})))=
H_{Q}(\Gamma(X,\Omega_{X}^{ch}))^{\ast}=H_{DR}(X)^{\ast}.
\eqno{(2.6.1)}
$$

{\bf 2.7.} One remark  concerning the ``size'' of
$D(\Gamma(X,\Omega_{X}^{ch})))$ is in order. This module also has
a gradation by conformal weight
$$
D(\Gamma(X,\Omega_{X}^{ch}))=\oplus_{n\geq 0}
D(\Gamma(X,\Omega_{X}^{ch}))_{n},\;
D(\Gamma(X,\Omega_{X}^{ch}))_{n}=
((\Gamma(X,\Omega_{X}^{ch})_{n})^{\ast}.
$$

One can argue that the space $D(\Gamma(X,\Omega_{X}^{ch}))_{n}$ is
often much
``bigger'' than  $\Gamma(X,\Omega_{X}^{ch})_{n}$. For example, since
$\Gamma(X,\Omega_{X}^{ch})_{0}$ is the de Rham complex of $X$, the dual
$D(\Gamma(X,\Omega_{X}^{ch}))_{0}$ is the space of distributions.
If $X$ is a compact manifold and we work in the $C^{\infty}$-setting,
then there is the following embedding
$$
\Gamma(X,\Omega_{X}^{ch})_{0}\hookrightarrow
D(\Gamma(X,\Omega_{X}^{ch}))_{0},\; \omega(\nu)=\int_{X}\omega\wedge\nu.
\eqno{(2.7.1)}
$$

If  $\Omega^{ch}_{X}$ were a coherent locally trivial sheaf, then
the same procedure applied to all conformal weights
spaces would give an embedding
$\Gamma(X,\Omega_{X}^{ch})\hookrightarrow
D(\Gamma(X,\Omega_{X}^{ch}))$. It is not, but $\Omega^{ch}_{X}$
does possess a filtration so that the associated graded sheaf is
indeed coherent and locally trivial. Therefore we can define a smaller
space $\tilde{D}(\Gamma(X,\Omega_{X}^{ch})))\subset
D(\Gamma(X,\Omega_{X}^{ch})))$ so that the graded object
$Gr( \tilde{D}(\Gamma(X,\Omega_{X}^{ch}))))$ with
 respect to the corresponding
dual filtration is the space of distributions of the type (2.7.1).
It is easy to see that $\tilde{D}(\Gamma(X,\Omega_{X}^{ch})))$ is
also a $\Gamma(X,\Omega_{X}^{ch})$-module.

\bigskip

\newpage

\centerline{\bf \S 3 Correlation functions}

\bigskip

In this section we shall keep to the vertex algebra
$\Gamma(X,\Omega_{X}^{ch})$ in the $C^{\infty}$-setting.

{\bf 3.1.}  Let $M_{1},..,M_{m}$ be restricted
$\Gamma(X,\Omega_{X}^{ch})$-modules. Attaching $M_{i}$ to
$P_{i}\in\CC$ we get a bundle over the configuration space
with infinite dimensional fiber $\otimes_{i}M_{i}$. Since
$z\BC[[z]]d/dz$ acts trivially on the conformal weight zero
component of each $M_{i}$, this bundle has a trivial subbundle
with fiber $\otimes_{i}(M_{i})_{0}$. Passing to
$\langle M_{1},...,M_{m}\rangle_{\CC}$, as explained in 1.6, we get a
trivial quotient bundle, to be denoted by
$\langle M_{1},...,M_{m}\rangle_{\CC,0}$, with fiber ``the space of
$\Lie( \Gamma(X,\Omega_{X}^{ch}))_{out}$-invariant functionals
restricted to  $\otimes_{i}(M_{i})_{0}$.

Define $Z(M)=Ker(Q_{0}:M\rightarrow M)$.
We can further quotient the bundle
$\langle M_{1},...,M_{m}\rangle(0)_{\CC}$ to get the bundle
$\langle Z(M_{1}),...,Z(M_{m})\rangle_{\CC,0}$ whose fiber
is a result of restricting that of
$\langle M_{1},...,M_{m}\rangle_{\CC,0}$ to
$\otimes_{i}Z(M_{i})_{0}$.

We are now in a position to reproduce a well-known physics
calculation that shows, in our terminology,
 that the connection on $\langle M_{1},...,M_{m}\rangle_{\CC}$
(see (1.6.4)) descends on the bundle
$\langle Z(M_{1}),...,Z(M_{m})\rangle_{\CC,0}$ and the
result is the trivial connection. Let $\langle.,.,...,.\rangle$ be
an element of the fiber of $\langle M_{1},...,M_{m}\rangle$;
apply  the vertical component of the connection $\Delta_{t}$
and restrict the result to $\otimes_{i}(Z(M_{i}))_{0})$. We have
$$
\langle \omega_{1},...,L_{-1}\omega_{t},\omega_{t+1},...,\omega_{m}
\rangle=
-\langle \omega_{1},...,[Q_{0},G_{-1}]\omega_{t},
\omega_{t+1},...,\omega_{m}
\rangle
$$
$$
=\Sigma_{s\neq t}\pm
\langle \omega_{1},...,Q_{0}\omega_{s},...,G_{-1}\omega_{t},
\omega_{t+1},...,\omega_{m}
\rangle=0
\eqno{(3.1.1)}
$$
Here we have used that , firstly, $L_{-1}=-[Q_{0},G_{-1}]$ by
(2.4.3b), secondly, that
 $Q_{0}\in\Lie(\Gamma(\Omega^{ch}_{X}))_{out}$ by (2.4.4) and,
finally, that $\omega_{i}\in Z(M_{i})$.

Similar calculation shows that
$$
\langle \omega_{1},...,\omega_{m}
\rangle= 0 \text{ if for some }i,\; \omega_{i}\in Im Q_{0} .
\eqno{(3.1.2)}
$$

{\bf 3.2.} Let us see what this construction gives us when
$\CC=\cp$ and the modules are either $\Gamma(X,\Omega^{ch}_{X})$ or
$D(\Gamma(X,\Omega^{ch}_{X}))$. We shall keep to the case when
$P_{1}=\infty$, $P_{i}=z_{i}, i=2,...,m-1$, $P_{m}=0$, $M_{1}=
D(\Gamma(X,\Omega^{ch}_{X}))$, $M_{i}=\Gamma(X,\Omega^{ch}_{X})$.
 In this case
the bundle
$\pi^{\ast}
\langle D(\Gamma(X,\Omega^{ch}_{X})),\Gamma(X,\Omega^{ch}_{X}),...,
\Gamma(X,\Omega^{ch}_{X})\rangle_{\CC}$
 is trivial (see (1.9.5) for the definition of $\pi$) because
our standard choice of the local coordinate $z-z_{i}$ at the point $z=z_{i}$
depends on the point ``smoothly''.
Hence any section of this bundle can be regarded as a
functional, depending on $z_{2},...,z_{m-1}$, on
$D(\Gamma(X,\Omega^{ch}_{X}))\otimes
\Gamma(X,\Omega^{ch}_{X})^{\otimes(m-1)}$.
 However the connection is not trivial
and the functional $\Phi_{z_{2},...,z_{m-1}}(...)$ defined by (1.9.1)
is a horizontal section of this bundle non-trivially dependng on
$z_{2},...,z_{m-1}$.
 What(3.1.1)
 tells us is that when restricted to the conformal
weight 0 components, the functional $\Phi_{z_{2},...,z_{m-1}}(...)$
is actually independent of $z_{2},...,z_{m-1}$. In fact one can avoid
using (3.1.1): given
$\omega_{1}^{\ast}\in D(\Gamma(X,\Omega^{ch}_{X}))_{0},
\omega_{i}\in\Gamma(X,\Omega^{ch}_{X}),\; i=2,...m$, it is easy to calculate
$\Phi_{z_{2},...,z_{m-1}}(\omega_{1}^{\ast},\omega_{2},...,\omega_{m})$
by directly using (1.9.1) and obtain:
$$
\Phi_{z_{2},...,z_{m-1}}(\omega_{1}^{\ast},\omega_{2},...,\omega_{m})
=
\omega_{1}^{\ast}(\omega_{2}\wedge...\wedge\omega_{m}).
\eqno{(3.2.1)}
$$
In r.h.s. of (3.2.1) the wedge product of forms is used; it makes
sense for, as we reminded the reader in 2.6, the
conformal weight 0 component of $\Gamma(X,\Omega^{ch}_{X})$ is exactly
the space of global differential forms.

Assuming further that $\omega_{1}^{\ast}$ actually belongs to
$\tilde{D}(\Gamma(X,\Omega^{ch}_{X}))_{0}$ (see 2.7), that is, equals
a  differential form, one rewrites (3.2.1) in the following nicer form:
$$
\Phi_{z_{2},...,z_{m-1}}(\omega_{1}^{\ast},\omega_{2},...,\omega_{m})
=\int_{X}
\omega_{1}^{\ast}\wedge\omega_{2}\wedge...\wedge\omega_{m}.
\eqno{(3.2.2)}
$$

The last formula suggests to regard the map
$$
\Phi_{\emptyset}:\tilde{D}(\Gamma(X,\Omega^{ch}_{X}))\otimes
\Gamma(X,\Omega^{ch}_{X})\rightarrow\BC,
\eqno{(3.2.3)}
$$
arising in the case of 2 points, as a ``chiralization'' of the standard
pairing of forms.

Similarly, in the case of 3 points, $\Phi_{z}(\omega_{1}^{\ast},
\omega_{2},\omega_{3})$ gives the structure constants of the algebra
of global differential forms.

Finally, by (3.1.2) and (2.6.1), $\Phi_{z_{2},...,z_{m-1}}(....)$ descends
 to a functional
$$
\Psi_{m}: H_{DR}(X)^{\ast}\otimes H_{DR}(X)^{\otimes (m-1)}\rightarrow
\BC.
\eqno{(3.2.4)}
$$

The latter is something very well known: $\Psi_{2}(.,.)$ is the Poincare
duality, $\Psi_{3}(,.,.,.)$ gives the structure constants of the cohomology
ring of $X$ etc.

\newpage

%


 \centerline{\bf Part III:  Calculation of
$\Gamma(\cpn,\Omega^{ch}_{\cpn})$}

\bigskip

\bigskip

\centerline{\bf \S 1. An embedding $L(Vect_{X})
\hookrightarrow \End(\Omega_{X}^{ch})$}

\bigskip

{\bf 1.1.} Let $X$ be a smooth manifold, $TX$ the total space of the
tangent bundle and $\Pi TX$ the supermanifold obtained by changing
the parity of all the fibers of the projection $ TX\rightarrow X$;
in other words, $\Pi TX$ is that supermanifold whose structure sheaf
is the sheaf of differential forms on $X$. The action of vector fields
on forms by the Lie derivative gives an embedding  of Lie
algebras:
$$
\pi: Vect_{X}\hookrightarrow Vect_{\Pi TX},
\eqno{(1.1.1)}
$$
which in coordinates reads as follows:
$$
\pi(f^{i}(b)a^{i})=f^{i}(b)a^{i}+\frac{\partial f^{j}(b)}{\partial b^{s}}
\phi^{s}\psi^{j}.
\eqno{(1.1.2)}
$$
In the last formula we keep to the following conventions:
$\{b^{i}\}$ is a coordinate system, $\phi^{i}=db^{i}$, so that
$\{b^{i},\phi^{i}\}$ is a coordinate system on $\Pi TX$; $a^{i}$
stands for the even vector field $\partial/\partial b^{i}$, while
$\psi^{i}$ signifies the odd one $\partial/\partial \phi^{i}$.

{\bf 1.2.}  Let us do a chiral version of 1.1.

  Let $Vect_{\BC^{N}}$ be the Lie algebra of polynomial vector fields
 on $\BC^{n}$
and let
$$
L(Vect_{\BC^{N}})=Vect_{\BC^{N}}\otimes\BC[t,t^{-1}]
\eqno{(1.2.1)}
$$
be the corresponding loop algebra. For $\tau\in Vect_{\BC^{N}}$ set
$\tau_{n}=\tau\otimes t^{n}\in L(Vect_{\BC^{N}})$.

Consider the vertex algebra $\Omega_{N}$, which
 was reviewed in II.2.4, and the space
of fields $\Fields(\Omega_{N})$ associated with $\Omega_{N}$, which
 was defined in I.1.4.
Attach to $\tau=f^{i}(b)a^{i}\in Vect_{\BC^{N}}$ a field
$\tau(z)\in\Fields(\Omega_{N})$:
$$
\tau(z)=:a^{i}(z)f^{i}(b(z)):+
:\frac{\partial f^{j}}{\partial b^{s}}(b(z))
\phi^{s}(z)\psi^{j}(z): .
$$
Let $\{\tau(z)_{n},n\in\BZ\}$ be Fourier coefficients of $\tau(z)$ so
that (in accordance with (I.1.6.2))
$$
\tau(z)=\Sigma_{n\in\BZ}\tau(z)_{n}z^{-n-1}.
$$

{\bf Lemma.} {\it The map
$$
\hat{\pi}: L(Vect_{\BC^{N}})\rightarrow \End(\Omega_{N}),
\tau_{n}\mapsto \tau(z)_{n}
\eqno{(1.2.2)}
$$
is a Lie algebra homomorphism.}

{\bf Proof.}    We have to check that
$$
[\hat{\pi}(\tau_{n}),\hat{\pi}(\xi_{m})]=\hat{\pi}([\tau,\xi]_{n+m}).
\eqno{(1.2.3)}
$$
In order to do so we compute the  OPE of the corresponding elements of
$\Fields(\Omega_{N})$,
$$
\tau(z)=:a^{i}(z)f^{i}(b(z)):+
:\frac{\partial f^{j}}{\partial b^{s}}(b(z))
\phi^{s}(z)\psi^{j}(z):
$$ and
 $$
\xi(z)=:a^{i}(z)g^{i}(b(w)):+
:\frac{\partial g^{j}}{\partial b^{s}}(b(z))
\phi^{s}(z)\psi^{j}(z):,
$$
by making use of Wick's theorem. In this way we get 2 singular terms:
one having pole of degree 2 and another one having pole of degree 1.
The degree 2 pole arises because of double pairings between $a$ and $b$ in
the first summands of each of the expressions above and because
of double pairings between $\phi$ and $\psi$ in the second summands;
as a result we have:
$$
\frac{1}{(z-w)^{2}}\{-\frac{\partial f^{i}}{\partial b^{l}}(b(z))
\frac{\partial f^{l}}{\partial b^{i}}(b(w))+
\frac{\partial f^{j}}{\partial b^{s}}(b(z))
\frac{\partial f^{s}}{\partial b^{j}}(b(w))\}=0.
$$

The degree one pole is caused by single pairings. The corresponding
calculation is parallel to the classical calculation of the bracket
$[\tau,\xi]$; here is the result:
 $$
\frac{[\tau,\xi](w)}{(z-w)}.
$$
Therefore
$$
\tau(z)\xi(w)\sim
\frac{[\tau,\xi](w)}{(z-w)}.
$$
Using (I.1.6.10) we get the desired equality (1.2.3).

{\bf 1.3.} {\bf Corollary.} {\it On a smooth manifold $X$, the homomorphisms
$\hat{\pi}$ glue to the homomorphism:
$$
\hat{\pi}_{X}:L(\Gamma(X, TX))\rightarrow \End(\Omega^{ch}_{X})
\eqno{(1.3.1)}
$$}

(See [MSV] or II.2.5 for the definition of the sheaf $\Omega^{ch}_{X}$.)

{\bf 1.4.} An attempt to replace $\Omega_{N}$
with a purely even vertex algebra $\CO_{N}$ (see II.2.1)
 in Lemma 1.2 fails:
as we showed in [MSV], a mapping similar to (1.2.2) may only give
an embedding of a certain extension of $L(Vect_{\BC^{N}})$.

\bigskip

\bigskip

\centerline{\bf \S 2 Calculation of $\Gamma(\cpn,\Omega_{\cpn}^{ch})$}

\bigskip

{\bf 2.1.} Let $G$ be a simple complex Lie group
and $\fg=\Lie G$.
It immediately follows from 1.3 that if $X$ is a $G$-space
or, at least, if there is a Lie algebra homomorphism:
$$
\rho:\fg\rightarrow \Gamma(X, TX),
\eqno{(2.1.1)}
$$
then $\Omega_{X}^{ch}$ is a sheaf of $\hfg$-modules of level 0 due to the
map
$$
\hat{\pi}_{X}\circ\hat{\rho}:\hfg\rightarrow\End(\Omega_{X}^{ch}),
\eqno{(2.1.2a)}
$$
where $\hat{\rho}$ is the obvious continuation of $\rho$ to the homomorphism
of loop algebras:
$$
\hat{\rho}:L(\fg)\rightarrow L(\Gamma(X, TX)),\;
\hat{\rho}(g\otimes t^{n})=\rho(g)\otimes t^{n}.
\eqno{(2.1.2b)}
$$

Focus on the case of $X=G/P$, where
 $P\subset G$
 is a parabolic subgroup. Let $X_{0}\in G/P$ be a big cell. Then
there is a $\hfg$-module embedding:
$$
\Gamma(X,\Omega_{X}^{ch})\hookrightarrow \Gamma(X_{0},\Omega_{X}^{ch}).
$$
One can think of $\Gamma(X,\Omega_{X}^{ch})$
as those sections of $\Omega_{X}^{ch}$ over $X_{0}$
that are ``regular outside
$X_{0}$''. A little thought shows that in fact
$$
\Gamma(X,\Omega_{X}^{ch})=\Gamma(X_{0},\Omega_{X}^{ch})^{\fg-\text{int}},
\eqno{(2.1.3)}
$$
where $\Gamma(X_{0},\Omega_{X}^{ch})^{\fg-\text{int}}$ stands for the
maximal $\fg$-integrable submodule or, equivalently, for the maximal
submodule on which $\fg$ operates locally finitely.

{\bf 2.2.} {\bf Example.} Let $X$ be $\cp$, $X_{0}=\BC$.
Let $\{E_{ij}\in gl_{2},\;1\leq i,j\leq 2\}$ be  the standard
basis of $gl_{2}$. In the homogeneous coordinates $(y_{0}:y_{1})$
 on $\cp$ one has the following specialization of (2.1.1):
$$
\rho: gl_{2}\rightarrow \Gamma(\cp,T\cp),
$$
$$\rho(E_{ij})=y_{i+1}\partial/\partial_{y_{j+1}}.
\eqno{(2.2.1)}
$$
Passing to the coordinate  $b=y_{1}/y_{0}$ on the big cell
$\BC\subset\cp$ and keeping to the conventions of 1.1, one has:
$$
\rho(E_{12})=a,\;\rho(E_{21})=-b^{2}a.
$$
Then (1.1.2) takes the following form
$$
\pi\circ\rho(E_{12})=a,\;\pi\circ\rho(E_{21})=-b^{2}a-2b\phi\psi.
\eqno{(2.2.2)}
$$
Passing further
to the loop algebras one obtains the following version of (2.1.2):
$$
\hat{\pi}\circ\hat{\rho}(E_{12}\otimes t^{n})=a(z)_{n},\;
\hat{\pi}\circ\hat{\rho}(E_{21}\otimes t^{n})=-:b(z)^{2}a(z):_{n}
-2:b(z)\phi(z)\psi(z):_{n}.
\eqno{(2.2.3)}
$$

Observe
that $J(z)=:\phi(z)\psi(z):$ is a ``free boson'' meaning that  the following
OPE is valid: $J(z)J(w)\sim J(w)/(z-w)^{2}$.
 The space $\Gamma(\BC,\Omega_{\cp}^{ch})$ is graded by fermionic
charge, that is to say, by eigenvalues of $J_{0}$:
 $\Gamma(\BC,\Omega_{\cp}^{ch})=
\oplus_{i\in\BZ}\Gamma(\BC,\Omega_{\cp}^{ch})^{(i)}$.
Hence, the boson-fermion
correspondence (see e.g. [K] 5.1) tells us that
$\Gamma(\BC,\Omega_{\cp}^{ch})^{(i)}$ is an irreducible representation
of the Lie algebra of Fourier components
 of the  fields
$a(z),b(z), J(z)$ generated by the vector $\phi_{-i+1}
\phi_{-i+2}\cdots\phi_{0}|0\rangle$ if $i\geq 0$, or
$\psi_{i}
\phi_{i-1}\cdots\phi_{-1}|0\rangle$ if $i<0$.

 With all of this in mind one checks (2.2.3) against (2.10) in [FF2] and
concludes that
 each $\Gamma(\BC,\Omega_{\cp}^{ch})^{(i)}$ is
a Wakimoto module over
$\widehat{sl}_{2}$ and the complex $(\Gamma(\BC,\Omega_{\cp}^{ch}),Q_{0})$
is the two-sided resolution of the trivial representation constructed
in [FF1]. Finally a glance at the diagram $III_{-}$
 in [FF2]4.2 allows one to use
(2.1.3) in order to obtain a rather explicit description
of the $\widehat{sl}_{2}$-module structure of
$\Gamma(\cp,\Omega_{\cp}^{ch}) $. Let us formulate
the result leaving the details of this
calculation out.

Let $V_{m}$ denote the simple $m+1$-dimensional $sl_{2}$-module.
Represent $\widehat{sl}_{2}$ as a direct sum of the loop algebra
$L(sl_{2})$ and $\BC\BK$, where $\BK$ is the standard central
element. The subalgebra $L(sl_{2})_{+}\subset \widehat{sl}_{2}$
of  loops regular at 0 maps onto $sl_{2}$ by means of the evaluation
at 0 map. Hence $V_{m}$ becomes an $L(sl_{2})_{+}$-module and one
defines the Weyl module $\BV_{m}$ of zero central charge as follows:
$$
\BV_{m}=\text{Ind}_{L(sl_{2})_{+}\oplus\BC\BK}^{\widehat{sl}_{2}} V_{m};
$$
it is assumed that $\BK$ operates on $V_{m}$ as 0.

$\BV_{m}$ has a unique irreducible quotient to be denoted $L_{m}$.
In fact, $\BV_{m}$ has a unique proper submodule isomorphic to
$L_{m+2}$.

{\bf Lemma.} {\it $\Gamma(\cp,\Omega_{\cp}^{ch})^{(i)}$ has a filtration by
$\widehat{sl}_{2}$-submodules
$$
F_{0}\subset F_{1}\subset\cdots,\; \cup_{m=0}^{\infty}F_{m}=
\Gamma(\cp,\Omega_{\cp}^{ch})^{(i)}
$$
so that:

if $i\leq 0$, then
$$
F_{0}=\BV_{-2i},\;\frac{F_{m}}{F_{m-1}}=\BV_{-2i+4m},m\geq 1;
$$
if $i> 0$, then
$$
F_{0}=L_{2i+2},\;\frac{F_{m}}{F_{m-1}}=\BV_{2i+4m},m\geq 1.
$$}

\bigskip

To further elaborate on
 the link between our approach and that of [FF2]
let us mention that the collection
of fields
$$
\{\psi(z),b(z)\psi(z),b(z)^{2}\psi(z)\}\subset
\End(\Gamma(X_{0},\Omega^{ch}_{\cp}))[[z,z^{-1}]]
\eqno{(2.2.4)}
$$
 is
the vertex operator associated with the adjoint representation
of $sl_{2}$ and that the
``chiral de Rham differential''
 $Q_{0}=\int a(z)\phi(z)$ coincides with the ``screening charge''.

In fact, the states $\psi_{-1}|0\rangle ,b_{0}\psi_{-1}|0\rangle,
b_{0}^{2}\psi_{-1}|0\rangle \in  \Gamma(X_{0},\Omega^{ch}_{\cp})$, to which
the fields $\psi(z),b(z)\psi(z),b(z)^{2}\psi(z)\}$ correspond
belong to $\Gamma(\cp,\Omega^{ch}_{\cp})$. Therefore (2.2.4)
 can be sharpened as follows
$$
\{\psi(z),b(z)\psi(z),b(z)^{2}\psi(z)\}\subset
\End(\Gamma(\cp,\Omega^{ch}_{\cp}))[[z,z^{-1}]].
\eqno{(2.2.5)}
$$

There is a well-known superaffine Lie algebra $(\widehat{sl}_{2})_{super}$
(of zero central charge), see e.g. [K] 2.5,
 obtained by taking the semi-direct
product of the usual loop algebra with the module of loops in the
 adjoint representation and declaring the
latter subspace odd. The formula (2.2.5) simply
means that the sheaf $\Omega^{ch}_{\cp}$ is actually a sheaf of
$(\widehat{sl}_{2})_{super}$-modules.

{\bf 2.3.} More generally one can consider $X$ equal $\cpn$,
$X_{0}=\BC^{N}$. By (2.1.2) $\Gamma(\cpn,\Omega_{\cpn}^{ch})$ is
a $\widehat{gl}_{N+1}$-module. Let $(y_{0}:...:y_{N})$ be the
homogeneous coordinates on $\cpn$; $b^{i}=y_{i}/y_{0}$, $i=1,...,N$
is a coordinate system on $X_{0}$. The formula (2.2.3) generalizes
as follows:
$$
\hat{\pi}\circ\hat{\rho}(E_{ij}\otimes t^{n})=
:b^{i-1}(z)a^{j-1}(z):_{n}+:\phi^{i-1}(z)\psi^{j-1}(z):_{n},\;
i,j\neq 1,
\eqno{(2.3.1a)}
$$
$$
\hat{\pi}\circ\hat{\rho}(E_{1j}\otimes t^{n})=
a^{j-1}(z)_{n},\;
j\neq 1,
\eqno{(2.3.1b)}
$$
$$
\hat{\pi}\circ\hat{\rho}(E_{i1}\otimes t^{n})= -
\Sigma_{l=1}^{N}
:b^{i-1}(z)b^{l}a^{l}(z):_{n}
$$
$$
-\Sigma_{l=1}^{N}:b^{i-1}(z)\phi^{l}(z)\psi^{l}(z):_{n}
-\Sigma_{l=1}^{N}:b^{l}(z)\phi^{i-1}(z)\psi^{l}(z):_{n}
,\;
i\neq 1.
\eqno{(2.3.1c)}
$$

 We again have $\Gamma(\BC^{N},\Omega_{\cpn}^{ch})=
\oplus_{i\in\BZ}\Gamma(\BC^{N},\Omega_{\cpn}^{ch})^{(i)}$. By
construction, each
\newline $\Gamma(\BC^{N},\Omega_{\cpn}^{ch})^{(i)}$   is
a generalized Wakimoto module attached in [FF1] to the maximal
parabolic subgroup of $SL_{N+1}$, or rather, (2.3.1a-c) provide
an explicit description of this module. By [MSV] Theorem 2.4 the complex
$(\Gamma(\BC^{N},\Omega_{\cpn}^{ch}),Q_{0})$ is a two-sided resolution
of the trivial representation composed of such modules.

By (2.1.3), $\Gamma(\cpn,\Omega_{\cpn}^{ch})$ is the maximal
$sl_{N+1}$-integrable submodule of the described generalized Waikmoto module.
To find an analogue of Lemma 2.2 one needs more information about generalized
Wakimoto modules than is avalaible now.

What was said in the end of 2.2 about vertex operators and the
structure of
a superaffine  algebra module on our sheaf carries over to the
present situation easily. For example, the vertex operator
associated with the adjoint representaion of $sl_{N+1}$
$$
\{e_{ij}(z), 1\leq i,j\leq N+1\}\subset
\End(\Gamma(\cpn,\Omega^{ch}_{\cpn}))[[z,z^{-1}]]
\eqno{(2.3.2)}
$$
is given (over $\BC^{N}$) by the following variation of (2.3.1):
$$
e_{ij}(z)=
:b^{i-1}(z)\psi^{j-1}(z):,\;
i,j\neq 1,
\eqno{(2.3.3a)}
$$
$$
e_{1j}(z)=
\psi^{j-1}(z),\;
j\neq 1,
\eqno{(2.3.3b)}
$$
$$
e_{i1}(z)= -
\Sigma_{l=1}^{N}
:b^{i-1}(z)b^{l}\psi^{l}(z):
\eqno{(2.3.3c)}
$$

\bigskip
\centerline{\bf References}
\bigskip

[B] R.~Borcherds, Vertex algebras, Kac-Moody algebras, and the Monster,
{\it Proc. Natl. Acad. Sci. USA}, {\bf 83} (1986), 3068-3071.

[BD] A.~Beilinson, V.~Drinfeld, Chiral algebras I, Preprint.


[BFM] A.~Beilinson, B.~Feigin, B.~Mazur, Introduction to algebraic
field theory on curves, Preprint.

[BS] A.~Beilinson, V.~Schechtman, Determinant bundles and Virasoro
algebras, {\it Comm. Math.Phys.} {\bf 118} (1988), 651-701.



[FF1] B.~Feigin, E.~Frenkel, Representations of affine Kac-Moody algebras
and bosonization, in: V.~Knizhnik Memorial Volume, L.~Brink, D.~Friedan,
A.M.~Polyakov (Eds.), 271-316, World Scientific, Singapore, 1990.

[FF2] B.~Feigin, E.~Frenkel, Affine Kac-Moody algebras and semi-infinite
flag manifolds, {\it Comm. Math. Phys.} {\bf 128} (1990), 161-189.





[K] V.~Kac, Vertex algebras for beginners, University Lecture Series,
{\bf 10}, American Mathematical Society, Providence, RI, 1997.





[MSV] F.~Malikov, V.~Schechtman, A.~Vaintrob, Chiral de Rham complex,
{\it Comm. Math. Phys.} (1999), to appear.




\bigskip

F.M.: Department of Mathematics, University of Southern California,
Los Angeles, CA 90089, USA;\ fmalikov\@mathj.usc.edu

V.S.: Department of Mathematics, Glasgow University,
15 University Gardens, Glasgow G12 8QW, UK;\
vs\@maths.gla.ac.uk

\enddocument